\documentclass{amsart}
\usepackage{amsmath, amsthm, amssymb, enumerate, verbatim,xcolor}
\usepackage[utf8]{inputenc}
\usepackage[english]{babel}

\usepackage{fixme}
\fxsetup{theme=color,mode=multiuser}
\definecolor{fxtarget}{rgb}{0.8000,0.0000,0.0000}

\newcommand{\andd}{\ensuremath{\qquad\mbox{and}\qquad}} 
\renewcommand{\epsilon}{\varepsilon}

\newcommand{\F}{\ensuremath{\mathbb{F}}}
\newcommand{\R}{\ensuremath{\mathbb{R}}}

\newcommand{\Z}{\ensuremath{\mathbb{Z}}}
\newcommand{\calL}{\ensuremath{\mathcal{L}}}
\newcommand{\calE}{\ensuremath{\mathcal{E}}}

\usepackage{bbold}
\newcommand\one{\mathbb{1}}

\theoremstyle{plain}
\newtheorem{thm}{Theorem}
\newtheorem{lem}[thm]{Lemma}
\newtheorem{prop}[thm]{Proposition}
\newtheorem{cor}[thm]{Corollary}
\newtheorem{claim}{Claim}
\newtheorem*{claim*}{Claim}

\theoremstyle{definition}

\newtheorem*{definition*}{Definition}

\theoremstyle{remark}

\renewcommand{\L}{\mathcal{L}}
\renewcommand{\P}{\mathcal{P}}

\newcommand{\bpm}{\begin{pmatrix}}
\newcommand{\epm}{\end{pmatrix}}

\newcommand{\SO}{\ensuremath{\mathrm{SO}_2(\F)}}
\newcommand{\SL}{\ensuremath{\mathrm{SL}_2(\F)}}
\newcommand{\SF}{\ensuremath{\mathrm{SF}_2(\F)}}
\newcommand{\GL}{\ensuremath{\mathrm{GL}}}
\newcommand{\B}{\ensuremath{\mathcal{B}}}
\newcommand{\I}{\ensuremath{\mathcal{I}}}
\newcommand{\Cc}{\ensuremath{\mathcal{C}}}
\newcommand{\PF}{\ensuremath{\mathbb{FP}^3}}

\newcommand{\PP}{\ensuremath{\mathbb{P}}}
\newcommand{\BG}{Blaschke-Gr\"unwald}
\newcommand{\T}{\ensuremath{\mathcal{T}^*}}

\newcommand{\bisector}{\ensuremath{\mathsf{B}}}
\newcommand{\cl}{\ensuremath{\mathcal{C}\ell}}

\newcommand{\fa}{\ensuremath{\mathfrak{a}}}
\newcommand{\fb}{\ensuremath{\mathfrak{b}}}
\newcommand{\fg}{\ensuremath{\mathfrak{g}}}
 
\newcommand{\fv}{\ensuremath{\mathfrak{v}}}

\def\Dpin{\Delta_{\text{pin}}}
\newcommand{\ax}{\ensuremath{\mathrm{Ax}}}

\title[Pinned Distances over Fields]{On the Pinned Distances Problem in Positive Characteristic}

\author{B. Murphy}
\address{Department of Mathematics, University of Bristol, UK}
\email{brendan.murphy@bristol.ac.uk}
\author{G. Petridis}
\address{Department of Mathematics, University of Georgia, Athens, GA 30602, USA}
\email{petridis@uga.edu}

\author{T. Pham}
\address{Department of Mathematics, University of Rochester, NY, USA}
\email{phamanhthang.vnu@gmail.com}

\author{M. Rudnev}
\address{Department of Mathematics, University of Bristol, UK}
\email{misha.rudnev@bristol.ac.uk}

\author{S. Stevens}
\address{Johann Radon Institute for Computational and Applied Mathematics (RICAM), Linz, Austria}
\email{sophie.stevens@ricam.oeaw.ac.at}

\date\today

\begin{document}
\begin{abstract} We study the Erd\H os-Falconer distance problem for a set $A\subset \F^2$, where  $\F$ is a field of positive characteristic $p$. If $\F=\F_p$ and the cardinality  $|A|$ exceeds $p^{5/4}$, we prove that $A$ determines an asymptotically full proportion of the feasible $p$ distances. 
For small sets $A$, namely when $|A|\leq p^{4/3}$ over any $\F$, we prove that either  $A$ determines $\gg|A|^{2/3}$ 
distances, or $A$ lies on an isotropic line. 
For both large and small sets, the results proved are in fact for pinned distances.  
\end{abstract}

\maketitle
\setcounter{tocdepth}{1}

\section{ Introduction}

For a compact set $A\subseteq \R^2$, the Falconer conjecture \cite{falconer86} claims that, for the distance set of $A$ to have positive Lebesgue measure, the Hausdorff dimension $d_H$ of $A$ must be strictly greater than $1$. 
Recently, Guth et al. \cite{guth2018falconer} made significant progress by  establishing the inequality $d_H\geq 5/4$ using decoupling.
This improved the estimate $d_H\geq 4/3$ due to Wolff \cite{wolff99falc} some 20 years previously. 

Iosevich and the fourth author \cite{iosevich2007distance} stated a finite field $\F_q$-variant of the Falconer conjecture. The Wolff exponent $4/3$ was subsequently proved by Chapman et al. \cite{chapman2012pinned} for $q\equiv 3 \mod(4)$; this constraint was removed by Bennett  et al. \cite{bennett2017group}. Here, we define the distance $d$ between two points $x=(x_1,x_2)$ and $y=(y_1,y_2)$ in the plane $\F_q^2$ to be
    \begin{equation} 
    	\label{eq:ddef}
        d\left(x,y\right):= (x-y)\cdot(x-y) = (x_1-y_1)^2+ (x_2-y_2)^2.
    \end{equation}
These approaches were based on Fourier analysis, with the underlying considerations similar to those used over the reals. In contrast, there is no positive characteristic analogue of the multiscale decoupling analysis, adapted by Guth et al. to the Falconer problem along the lines of the general method of Bourgain and Demeter \cite{bourgain2015proof}.  Moreover, the first and the second authors \cite{MP_Falc} showed that, when $|A| \leq q^{4/3}$, subspaces over subfields in $\F_q$ generally preclude one from having more than $q/2$ distinct distances.

Theorem \ref{thm:Pinned Falconer} in this paper claims the exponent $5/4$ in the context of the prime residue field $\F_p$. To achieve this, we use state-of-the art geometric incidence theory, combined with several new geometric set-up ideas. These include the classical {\em Gr\"unwald-Blaschke kinematic map} which provides the ``correct coordinate system'', in which one can efficiently estimate the key quantity of {\em bisector energy}. Both the Gr\"unwald-Blaschke kinematic map and the bisector energy enter the proof of the general field Theorem  \ref{thm: Pinned F}, which we present first. This theorem also establishes the exponent $2/3$ for the Erd{\H{o}}s distinct distance problem in $\F^2$. Although this exponent is weaker than its real counterpart, it represents the best known result in general fields.  

\medskip
The Falconer conjecture is a continuous version of the 1946 question of Erd{\H{o}}s \cite{erdos1946distances}, who conjectured that any $N$-point set $A$ in the real plane should determine at least $c\frac{N}{\sqrt{\log{N}}}$ distinct distances,  for some universal constant $c>0$. The square grid, which achieves this bound, is thus expected to be a minimiser. In 2010 Guth and Katz \cite{guth2015erdos}  essentially resolved the question, and proved that $N$ points in $\R^2$ determine at least $c\frac{N}{\log{N}}$ distances. 
See \cite{iosevichetalbook} for a fuller account and history of the problem prior to Guth and Katz's result.
 
The distinct distance problem can also be posed over arbitrary fields $\F$ with distance defined as in \eqref{eq:ddef}. We assume throughout that $\F$ has a positive characteristic $p$, although the results apply to zero characteristic by just omitting the constraints involving $p$. The characteristic $p$ serves as an asymptotic parameter, and our results are trivial for small $p$, in particular the special case $p=2$.
	    
In this notation, the distinct distance problem is to find a lower bound on the cardinality 
\[\Delta(A):=|\{d(a,b):\, a,b\in A\}|
\]
of the set of distances determined by a finite point set $A\subset \F^2$, with $|\cdot|$ denoting the cardinality of finite sets. The lack of order in $\F$ seriously limits the availability and applicability of the tools used and developed for its study over $\R$.

We use the usual notations $X \ll Y$, equivalently  $X=O(Y)$ if there exists  some absolute constant $c>0$ such that $X \leq cY$. Conversely,  $X \gg Y$,  equivalently  $X=\Omega(Y)$ means $Y\ll X$. The constant implicit in this notation may freely change from line to line.  
	
	Erd{\H{o}}s formulated other versions of the distinct distance problem, including the stronger \emph{pinned distance variant}. Given a set $A\subseteq \F^2$, and a point $x\in \F^2$ (`a pin'), the cardinality of the set of pinned distances at $x$ with respect to $A$ is
	\begin{equation} \label{eq:Dax}
	    \Delta(A;x):= |\{d(a,x)\colon a\in A\}|\,.
	\end{equation}
    The pinned distance count of $A$ is 
    \[ 
      \Dpin(A):=\max_{x\in A}\Delta(A;x)\,,
    \]
    and the pinned distance problem is to determine $\min_{|A|=N}\Dpin(A)$.  It is clear that $\Delta(A) \geq \Dpin(A)$.
    Erd\H{o}s \cite{erdos1975some} conjectured that asymptotically, the behaviour of $\Delta(A)$ and $\Dpin(A)$ should be the same, at least in the reals, i.e. $\Dpin(A)\gg \frac{|A|}{\sqrt{\log{|A|}}}$.

The pinned distance conjecture is still wide open even over the reals,  with the strongest bound $\Dpin(A)\gg  |A|^{.8641\ldots}$ due to Katz and Tardos \cite{katz-tardos} in the early 2000s. 
	Questions about distances have also been asked in $\R^2$ (and higher dimensions) for non-Euclidean strictly convex norms, as well as hyperbolic, spherical distances, etc. See e.g.
 \cite{valtr2005strictly}, \cite{matouvsek2011number}, \cite{roche-newton2012minkowski},   \cite{rudnev2016klein}. For more relatives of the distinct distance problem over the reals we recommend the survey of Sheffer \cite{sheffer2014distinct}. 
 
 Erd\H{o}s' combinatorial pinned distance argument in \cite{erdos1946distances} works equally well for non-zero distances in $\F^2$: his argument shows that every $A\subset \F^2$ determines either only the zero distance or at least $\gg |A|^{1/2}$ non-zero distances. (Throughout the rest of this section we assume by default that  $\Delta(A)\neq\{0\}$, that is $A$ does not lie on a single {\em isotropic line} in $\F^2$.) The first improvement of this  over the prime field $\F_p$ was obtained by  Bourgain, Katz and Tao \cite{bourgain2005sum-product}, who proved a non-quantitative non-trivial bound on $\Dpin(A)$, based on a non-trivial Szemer\'edi-Trotter type point-line incidence bound, which in turn followed from a sum-product estimate.
        
 
 A quantitative lower bound  $\Dpin(A)\gg |A|^{8/15}$ (for general $\F$, henceforth under suitable smallness conditions on $|A|$ in terms of $p$)
 was established by the fifth author and de Zeeuw  \cite{stevens2017improved}  by using their novel quantitative  Szemer\'edi-Trotter type point-line incidence theorem over a general $\F$, proved in the same paper. In the special case $A=X\times X$, the second author \cite{petridis2016pinned} proved a stronger result $\Dpin(A)\gg |A|^{3/4}.$
In the context of $\F=\F_p$, for $p\equiv 3 \bmod(4)$,  Iosevich, Koh and the third author \cite{iosevich2019perspective}  improved the $\Dpin(A)$ exponent to $\frac{1}{2} + \frac{69}{1558} = 0.5442\dots$.
    
The paper by Lund and the second author~\cite{lund2018pinned} proved a slightly weaker but more general result $\Dpin(A)\gg |A|^{20/37},$ for $A\subset \F^2$.    
The result was based on upper bounds for the  {\em bisector energy} quantity, introduced by Lund, Sheffer and de Zeeuw \cite{lund2016bis},  
that is, the number of pairs of segments with vertices in $A$ that are symmetric relative to some line, namely the bisector line. Using the bisector energy enabled Hanson, Lund and Roche-Newton \cite{hanson2016distinct} to extend the earlier Falconer conjecture threshold exponent $4/3$  by Chapman et al. \cite{chapman2012pinned}, and to improve the result of Bennett  et al. \cite{bennett2017group} to pinned distances over $\F_q$.

In this paper, we will show that these bisector energy bounds respond most naturally to the current state of the art of incidence tools over a general field $\F$. Note that these tools are much weaker than the incidence tools developed over the reals. Geometric incidence bounds (for sufficiently small sets in positive characteristic) over a general $\F$ are largely confined to the fourth author's point-plane theorem in $\F^3$ \cite{rudnev2018point-plane} or its descendants \cite{yazici2015growth,murphy2017second,stevens2017improved}.

One of the key observations in this paper is that the bisector energy of $A\subseteq \F^2$ can be estimated directly via a point-plane incidence, if one assumes a judicious geometric viewpoint.  This viewpoint is provided by the Blaschke-Gr\"unwald kinematic mapping, described in Section~\ref{sec:blaschegrunwald}. The mapping provides a natural embedding of the plane Euclidean motion group in the projective three-space. This fact was partially rediscovered by Elekes and Sharir \cite{elekes2011incidences} followed by Guth and Katz \cite{guth2015erdos}; it  played a major role in their resolution of the distinct distance problem in $\R^2$.

This paper proves the new bound $\Dpin(A)\gg |A|^{2/3}$ over a general $\F$, thus matching Moser's exponent \cite{Moser1952} over the reals from almost 70 years ago.

    

We also establish a new Falconer-type threshold in the special case $\F = \F_p$: we show that it is sufficient for a set $A\subseteq \F_p$ to satisfy  $|A|\geq p^{5/4}$ to yield a positive 
proportion of  possible distances. We restate that the latter exponent $5/4$ cannot hold over $\F_q$ in terms of $q$, where the exponent $4/3$ (from  \cite{chapman2012pinned}, \cite{bennett2017group}, \cite{hanson2016distinct}) is generally the best one possible, owing to the series of examples in \cite{MP_Falc}. Our numerology matches that of the recent paper by Guth et al. \cite{guth2018falconer} that broke the ice in the progress towards the Falconer problem in $\R^2$; Guth et al. improved the critical Hausdorff dimension $d_H$ from  $4/3$  to $5/4$.
 

    \medskip


\section{Results and discussion of techniques} 

\subsection{Main results}\label{sec:mainresults}
The main results of this paper are a consequence of a new bound on the number of isosceles triangles with vertices in $A$ for a set $A\subseteq \F^2$. This is the content of the forthcoming Theorem~\ref{thm:Isosc Triang}, whose statement we defer until the relevant quantities have been developed. The relation between the number of isosceles triangles and pinned distances is formalised by the Cauchy-Schwarz inequality, together with some technical considerations regarding the number of collinear points in $A$. 

Our bound on the number of isosceles triangles consists of three terms, with each term optimised for a particular situation. For clarity, we present the pinned distance result that corresponds to each term separately, since their ranges of applicability are quite distinct. 

Our first result applies to the large set case, that is the pinned version of the Erd\H{o}s--Falconer distance problem for $A\subseteq \F_p^2$: recall that the Erd\H{o}s--Falconer distance problem asks for the threshold of $|A|$ so that $A$ determines a positive proportion of all feasible $p$ distances.  

\begin{thm}[Pinned Erd\H{o}s--Falconer over $\F_p$]
\label{thm:Pinned Falconer}
    There exists an absolute constant $c$ with the following property. For all primes $p$ and all $A \subseteq \F_p^2$ such that $|A| \geq p^{5/4}$, we have 
    \[\Dpin(A) > cp\,.\] 
    
    Moreover, for all $0 < \varepsilon <1$ and all increasing unbounded functions $\omega: \R^+ \to \R^+$, there exists $p_0=p_0(\varepsilon, \omega)$ such that if $p>p_0$ is a prime and $A \subseteq \F_p^2$ satisfies $|A| > \omega(p) p^{5/4}$, then for at least $(1 - \varepsilon) |A|$ many $a \in A$ we have $\Delta(A;a) > (1 - 2\varepsilon) p$. 
\end{thm}

Theorem~\ref{thm:Pinned Falconer} implies, with a somewhat informal notation, that 
\[
    \lim_{p \to \infty} \frac{|A|}{p^{5/4}} = \infty \implies \lim_{p \to \infty} \frac{\Dpin(A)}{p} = 1.
\]
This improves the $|A| \geq p^{4/3}$ threshold obtained by  Hanson, Lund, and Roche-Newton \cite{hanson2016distinct}. The unpinned result was proven in \cite{bennett2017group, chapman2012pinned} and like \cite{hanson2016distinct} also holds in a finite field $\F_q$ (of characteristic $p$), with $q$ replacing $p$ in the bounds. 


Our second statement concerns sets $A\subseteq \F_p^2$ of intermediate cardinality $p\ll |A|<p^{5/4}$. This bounds the  number of pinned distances in terms of $|A|$ and $p$:

\begin{thm}[Pinned distances over $\F_p$]\label{thm: Pinned Fp only}
    Let $A\subset \F_p^2$ be a set of points satisfying $4p <|A|<p^{5/4}$.
    Then
    \begin{equation}\label{eq:dpin Fp only}
        \Dpin(A)\gg |A|^{4/3}p^{-2/3}\,.        
    \end{equation}
\end{thm}

Theorem~\ref{thm:Pinned Falconer} and Theorem~\ref{thm: Pinned Fp only} share a single  proof, delivering a new bound on the number of isosceles triangles in Theorem ~\ref{thm:Isosc Triang} below.

Finally, we address the case of small sets, where $p$ appears only as the applicability constraint. This accounts for the special case $\F=\F_p$, where one cannot have more then $p$ distinct distances. 
\begin{thm}[Pinned distances over $\F$]\label{thm: Pinned F}
    Let $A\subset \F^2$ be a set of points.
    If $\F$ has positive characteristic $p>0$, assume in addition that $|A|\leq p^{4/3}$. 
    Then either $A$ is contained in  a single isotropic line, so all pair-wise distances are zero, or
    \begin{equation}\label{eq:dpin F}
        \Dpin(A)\gg |A|^{2/3}\,.        
    \end{equation}
\end{thm}

In the context of $\F_p$,  the bound \eqref{eq:dpin F} is stronger than \eqref{eq:dpin Fp only} whenever $ |A|<p$. Note that in Theorems~\ref{thm:Pinned Falconer}~and~\ref{thm: Pinned Fp only} it is impossible for $A$ to be entirely contained in an isotropic line. A result similar to Theorem~\ref{thm: Pinned F} holds for different distances. Given a quadratic  form $f: \F^2 \to \F$ we define $d_f(x,y) = f(x-y)$. If $\F$ is algebraically closed, then $f$ is equivalent to $x_1^2+x_2^2$ and so Theorem~\ref{thm: Pinned F} applies to $d_f$.

\subsection{Discussion of techniques}
Our results rely on a new bound on the number of non-degenerate isosceles triangles with vertices lying in the point set of interest $A$. We describe the techniques proving Theorem~\ref{thm: Pinned F}, and then indicate the differences in method we use to prove Theorems~\ref{thm:Pinned Falconer}~and~\ref{thm: Pinned Fp only}.

Let $\T(A)$ be the number of \emph{non-degenerate isosceles triangles} with vertices in $A$. That is, 
\[\T(A) = |\{ (a,b,c)\in A^3 \colon  d(a,b)=d(a,c) \text{ and }(b-c)\cdot (b-c)\neq 0\}|\,.
\]
The reasons behind this choice are discussed in Section~\ref{sec:bisect-energy-isosc}.

\begin{thm}\label{thm:Isosc Triang}
    Let $p$ be a prime and $A \subseteq \F_p^2$. If $|A| \leq p^{4/3}$, then
    \begin{equation}\label{eq:T*bound}
        \T(A) - \frac{|A|^3}{p} \ll \min\left(p^{2/3} |A|^{5/3} + p^{1/4} |A|^2, |A|^{7/3}\right)\,.
    \end{equation}
    Moreover, the bound $\T(A)\ll|A|^{7/3}$ holds for $A\subseteq \F^2$ whenever $|A|< \sf{char}(\F)^{4/3}$.
\end{thm}

We choose to write the bound \eqref{eq:T*bound} as an asymptotic formula. The $|A|^3/p$ term dominates the right-hand side for $|A|\geq p^{5/4}$. On the right-hand side, the $p^{2/3} |A|^{5/3}$ dominates when $p <  |A| \leq  p^{5/4}$ and the $|A|^{7/3}$ term dominates when $|A| \leq p$. We include the $p^{1/4} |A|^2$ term because it dominates the right-hand side when $|A| > p^{5/4}$. As an estimate, the bound \eqref{eq:T*bound} could be more simply restated as
\[
\T(A) \ll \left\{ \begin{array}{lll} |A|^3 / p &,& \text{if } |A| > p^{5/4} \\ p^{2/3} |A|^{5/3} &,& \text{if } p \leq |A| \leq p^{5/4} \\ |A|^{7/3} &,& \text{if } |A| < p \end{array} \right., 
\]
for all $A \subseteq \F_p^2$ (see \cite{hanson2016distinct} for the case $|A| > p^{4/3}$).

 The $|A|^3/p$ term is the expected number of ordered (non-degenerate) isosceles triangles if we choose $A$ randomly by including each element of $\F_p^2$ independently with probability $|A| / p^2$. 
Theorem~\ref{thm:Isosc Triang} can therefore be interpreted as saying that all sufficiently large sets in $\F_p^2$  exhibit pseudorandom behaviour with respect to $\T(A)$. 
 
The two minimands appearing in Theorem~\ref{thm:Isosc Triang} arise from the same techniques: we  follow the established method of studying isosceles triangles through the perpendicular bisectors of pairs of points in $A$.

The perpendicular bisector of points $a,b\in \F^2$ with $d(a,b)\neq 0$ is the line
    \begin{align*}
        \bisector(a,b) & =\{x\in \F^2 \colon d(a,x) = d(b,x)\} \\ &= \{x\in \F^2 \colon 2 x \cdot (a-b) = d(a,0) - d(b,0)\}.
    \end{align*}
We relate the number of bisectors determined by $A$ to (a subtle variant of) the \emph{bisector energy}\/ of the set $A$, which is the number of pairs of points whose perpendicular bisectors coincide:
    \[
        \B(A):=|\{(a,b,c,d)\in A^4\colon\bisector(a,b) = \bisector(c,d)\}|.
    \]
    In  simple terms, $\B(A)$ is the number of pairs of segments with endpoints in $A$, with the two segments symmetric relative to a reflection in some axis, namely the bisector. It is easy to see (and is shown below) that the axis cannot be isotropic, else a symmetry is not well-defined. 
    
The variant of the bisector energy that we use in the sequel, denoted as $\B^*(A)$, see Section~\ref{sec:bisect-energy-isosc}, imposes an additional counting restriction that pairs of segments, symmetric relative to the bisector, and thus having the same  length, are themselves non-isotropic. Quite clearly, this will not affect $\B(A)$ in a significant way.

The bisector energy controls the number of \emph{isosceles triangles}\/ in $A$; upper bounds for the number of isosceles triangles in $A$ yield lower bounds for $\Dpin(A)$.

Lund and the second author showed quantitatively that if the bisector energy -- in $\R^2$ -- is large, then $A$ must contain many collinear points or many co-circular points \cite[Theorem 1.2]{lund2018pinned}. Earlier Lund, Sheffer and de Zeeuw \cite[Theorem 2.2]{lund2016bis} gave an example of a (real plane)  point set with large bisector energy, the set being supported on a family of rich in points parallel lines. (One can modify the construction to turn parallel lines into concentric circles.) Here we succeed in showing that the bound of \cite[Theorem 1.2]{lund2018pinned} extends, in fact, to all fields $\F$. Furthermore,  in the context of $\F=\F_p$ and the cardinality regime pertaining to the Falconer conjecture, we also show a structural result: the only scenario that can account for large bisector energy is when a large proportion of $A$ is supported on families of rich parallel lines or concentric circles. It is easy to see that this structure contributes acceptably few  isosceles triangles.

We count the bisector energy by studying reflective symmetries: a quadruple in the support of $\B^*(A)$ is a pair of segments related by a reflection over a perpendicular bisector. The two segments in this support are of the same length.
We partition pairs of segments according to their length.
To prove our bound on the modified bisector energy $\B^*(A)$, we use the \emph{kinematic mapping}\/ of Blaschke and Gr\"unwald \cite{blaschke1960kinematik, grunwald1911abbildungsprinzip} to embed the space of segments of the same non-zero length into projective three-space.
The bisector energy in a class of $n$ segments of the same non-zero length is then represented by the number of incidences between $n$ points and $n$ planes, which we bound using the point-plane incidence theorem of the fourth author.\footnote{It is interesting to note that this reasoning is implicit in Lund and Petridis' proof of \cite[Theorem 1.2]{lund2018pinned}, which is developed in $\R^2$ and uses inversion in a circle instead of the Blaschke-Gr\"unwald kinematic mapping and Beck's theorem instead of Theorem  \ref{t:rudnev}. The proof also shows implicitly how the Beck theorem, which unfortunately is so far unavailable over the general $\F$ would easily imply Theorem \ref{t:rudnev} in $\R^3$, see \cite[Appendix]{petbis19}.}
    
To be precise, we use $S_r=S_r(A)$ to denote the set of pairs of points of distance $r$ apart:
    \begin{equation} \label{eq: Sr}
      S_r = S_r(A) :=\{(a,b)\in A^2\colon d(a,b)=r, a \neq b\}.
    \end{equation}
The modified bisector energy $\B^*(A)$ is equal to the sum over $r\not=0$ of the number of pairs of segments in $S_r(A)$ that are \emph{axially symmetric}\/ (plus an error term for isotropic segments).
As mentioned above, we count the number of such pairs by representing it as a point-plane incidence count in projective three-space; see Claim 1 below.
From this it follows that
    \[
        \B^*(A)\ll \sum_{r\not=0}|S_r|^{3/2} \ll  \left( \sum_{r\not=0}|S_r| \right)^{1/2} \left( \sum_{r\not=0}|S_r|^2 \right)^{1/2} \leq |A| \left( \sum_{r\not=0}|S_r|^2 \right)^{1/2},
    \]
unless the mapping of $S_r(A)$ into $\PF$ has too many collinear points. In terms of the original set $A$, this roughly translates into $A$ having many collinear or co-circular points.
In the part of the proof exploiting the kinematic mapping, we assume without loss of generality that $\F$ is algebraically closed; we may embed $A$ into the algebraic closure of $\F$ without decreasing the quantities we wish to bound. This idea proves the second minimand in Theorem~\ref{thm:Isosc Triang}, leading to Theorem~\ref{thm: Pinned F}.

In the proof of the first minimand in Theorem~\ref{thm:Isosc Triang}, and thus Theorems~\ref{thm:Pinned Falconer}~and~\ref{thm: Pinned Fp only}, we return to the ``configuration space'' $\F_p^2$. We observe that when proving the second minimand, we are in a bad situation if $A$ has many collinear or co-circular points. More precisely, this bad situation occurs when $A$ contains many segments of the same length with endpoints lying on pairs of parallel lines or concentric circles. Then, by the trivial fact in $\F = \F_p$ that the number of collinear or co-circular points is necessarily bounded in terms of $p$ (in general, by $p+1$), we obtain   a sufficient numerical advantage in this bad setting, once the isosceles triangles to be counted have been properly ``pre-pruned''.

To this effect, we decompose the count of isosceles triangles: we first separate out the count of isosceles triangles into those with an axis of symmetry which interacts (in a particular way to be later described) with rich lines or circles -- that is, lines or circles containing many points of $A$. To these triangles, we apply a counting argument to show that, since there cannot be too many rich lines or circles, the contribution of these types of triangles to $\T(A)$ is controlled. The count of the remaining contribution to $\T(A)$ proceeds again by a bisector argument. However, this time we are in a better situation to apply the point-plane incidence theorem of the fourth author (the version of restricted incidences \cite[Theorem 3*]{rudnev2014number}).

    
    
\section{Preliminaries}
\subsection{Distance preserving transformations}
    Let $\SO\subseteq\SL$ denote the set of unit determinant linear transformations preserving the distance:
    \[
        \SO:=\{g\in \SL\colon \forall x,y\in\F^2,\,d(x,y)=d(gx,gy)\}.
    \]
    As a matrix group,
    \[
        \SO  =\left\{
        \begin{pmatrix}
            u & - v \\
            v & u\\
        \end{pmatrix}:
        u,v\in \F, u^2+ v^2=1 \right\}.
    \]
    We will use the notation $\Cc\subseteq \F^2$ for the unit circle, and write $(u,v)\in \Cc$.
    As is the case for rotations acting on circles in $\R^2$, the group $\SO$ acts simply transitively on the level sets $\{(x,y)\in\F^2\colon d(x,y)=t\}$ for all $t\not=0$.
    Thus $d(x,y)=d(x',y') \neq 0$ if and only if there is a rotation $\theta\in \SO$ such that $\theta x - \theta y = x'-y'$.
    
    Let $T_2(\F)$ be the group of translations $x\mapsto x+t$ acting on the plane $\F^2$. 
    The group $\SF$ of positively oriented rigid motions of $\F^2$ is generated by $\SO$ and $T_2(\F)$; this is the analogue of the special Euclidean group $\mathrm{SE}_2(\R)$. 

It is well known that there is an injective group homomorphism from $\SF$ (where the group operation is composition of maps) into $\mathrm{SL}_3(\F)$ (where the group operation is matrix multiplication). Thus, an element of $\SF$ can be represented as a matrix of the form: 
    \begin{equation}\label{matrix}
      \begin{pmatrix}
        u&- v&s\\
        v&u&t\\
        0&0&1
           \end{pmatrix},
       \quad \mbox{where $u^2+ v^2 = 1$}.
    \end{equation}
    
    By the above discussion, we see that $d(x,y)=d(x',y')$ if and only if there exists $g\in \SF$ such that $g(x,y)=(x',y')$.
    If such a $g$ exists, an easy calculation shows that it is unique.

\subsection{Blaschke-Gr\"unwald kinematic mapping}\label{sec:blaschegrunwald}
    The Blaschke-Gr\"unwald kinematic mapping \cite{blaschke1960kinematik,grunwald1911abbildungsprinzip} assigns to each element
    $g\in \mathrm{SE}_2(\R)$ 
    a point in projective space $\R\mathbb{P}^3$.
    For a detailed exposition concerning this mapping and its properties, see the textbook by Bottema and Roth \cite[Chapter 11] {bottema1990theoretical}.
    The kinematic mapping was partially rediscovered some 100 years later by Elekes and Sharir \cite{elekes2011incidences} and played an essential role in the resolution of the Erd\H{o}s distinct distance problem in $\R^2$ by Guth and Katz \cite{guth2015erdos}.
    
    The definition of the original \BG{} kinematic mapping extends to all fields that are closed under taking square roots.
    The reason for this is the necessity to have well-defined  ``half-angles'': for all $(u,v)\in \mathcal C$ (the unit circle), we may resolve the system of quadratic equations
    \begin{equation}\label{trig}
        u = \tilde u^2 - \tilde v^2,\qquad v = 2\tilde u \tilde v\,.
    \end{equation}
    We choose a root of the equation $\tilde u^2 = \frac{1+u}{2}$ as a solution; since we use projective coordinates, it does not matter which of the two roots one chooses for $\tilde u$ -- this choice, once made, defines $\tilde v$ unambiguously. It then follows that $(\tilde u, \tilde v)\in \mathcal C$.
    With these preliminaries in hand, we may define the Blaschke-Gr\"unwald kinematic mapping, which embeds $\SF$ into $\PF$: an element of $\SF$ of the form of \eqref{matrix} becomes the projective point:
    \begin{equation}\label{bgr}
        [X_0:X_1:X_2:X_3] = [2 \tilde u: 2 \tilde v: s  \tilde u  +  t \tilde v \colon s \tilde v - t  \tilde u].
    \end{equation}
    Note that the mapping  \eqref{bgr} does not depend on the sign choice in the half-angle formulae \eqref{trig}.
    Conversely,
    \begin{equation}\label{bgrinv}
        u = \frac{X_0^2-X_1^2}{X_0^2+X_1^2},\; v=  \frac{2X_0X_1}{X_0^2+X_1^2} ,\; \frac{s}{2} = \frac{X_1X_3+X_0X_2}{X_0^2+X_1^2}, \; \frac{t}{2}=\frac{X_1X_2-X_0X_3}{X_0^2+X_1^2}.
    \end{equation}
    
    If $\F$ is a field where some elements do not have square roots, we can use projectivity to avoid them.
    If $\tilde u\not=0$, we may multiply the coordinates of the left hand side of \eqref{bgr} to find
    \[
        [X_0:X_1:X_2:X_3] = [2(u+1) : 2v : s(u+1)+  tv : s v - t(u+1)].
    \]

    Observe that the image of the kinematic mapping $\kappa$, is $\mathbb{FP}^3\setminus \{X_0^2+X_1^2=0\}$. That is one removes from $\mathbb{FP}^3$ the exceptional set, which is a line if $-1$ is not a square in $\F_p$ and is a union of two planes if $-1$ is a square in $\F_p$.
    
    The kinematic mapping has a number of remarkable properties, however, the easiest way to derive these properties is by studying a certain Clifford algebra.
    Since we do not have a reference for these computations over arbitrary fields, we provide them in Appendix~\ref{sec:clifford-calc}.
    
    The most important property of $\kappa$ for this paper is that translation in the group $\SF$ corresponds to a \emph{projective transformation}\/ of \PF.
    \begin{prop}
      \label{prop:kinematic-equivariant}
      For all $g\in\SF$ there are projective maps $\phi_g\colon \PF\to \PF$ and $\phi^g\colon \PF\to \PF$ such that for all $x\in\SF$
      \[
    \kappa(gx) =\phi_g(\kappa(x))\andd \kappa(xg)=\phi^g(\kappa(x)).
        \]
            \end{prop}
    The proof of this proposition is contained in the proof of Corollary~\ref{cor:1} in Appendix~\ref{sec:clifford-calc}.
    
    As a corollary, we see that the set of all rigid motions mapping one fixed point to another fixed point corresponds to a projective line.
    For points $x$ and $y$ in $\F^2$, let $T_{xy}$ denote the set of $g\in\SF$ such that $gx=y$.
    \begin{cor} \label{cor:transporters}
        For all $x$ and $y$ in $\F^2$, the image $\kappa(T_{xy})$ is a projective line.
    \end{cor}
    \begin{proof}
        The image of the rotation subgroup $\SO$ under $\kappa$ is $X_2=X_3=0$, which is a projective line.
        By the transformation properties, all conjugate subgroups of $\SO$ are projective lines, and all cosets of these groups are lines.
        The set $T_{xy}$ is a left coset of the stabiliser of $x$, which is conjugate to $\SO$.
    \end{proof}

\subsection{Isotropic lines} \label{sec:isotropy}
    In arbitrary fields, there may exist a set of points whose pairwise distance is 0. This is an obvious obstruction to obtaining a lower bound on $\Dpin(A)$, and so we have to consider these points separately. 
    
    A vector $v\neq 0\in \F^2$ is \emph{isotropic} if $d(v,0)=0$. If ${\rm i}:=\sqrt{-1}\in \F$, then $\F$ contains isotropic vectors. In particular, we note that when $p\equiv 3\mod(4)$, then $-1$ is not a square in $\F_p$ so there are no isotropic vectors. 
    Given a finite point set $A$, an oriented segment is a pair $(a,a')\in A^2$ with length $d(a,a')$. If $d(a,a')=0$, the segment is called isotropic; when $a\neq a'$ we say that $(a,a')$ is a non-trivial isotropic segment. A non-trivial isotropic segment lies on an isotropic line with slope $\pm {\rm i}$.
    
    Isotropic line segments should be excluded from counts, for there may be too many of them: a single isotropic line supporting $N$ points contains $\gg N^2$ zero-length segments. 
    
    Among other facts on isotropic lines, it is easy to see that if $(a,b,c)$ is an isosceles triangle with $d(a,b) = d(a,c) \neq 0$, then the perpendicular bisector $\bisector(b,c)$ is not isotropic \cite[Corollary 2.5]{lund2018pinned}.

\subsection{Axial symmetries}\label{sec:axial}
    As in the Euclidean case, $\SF$ has index two in the group of all distance-preserving transformations.
    The other coset of $\SF$ consists of compositions of reflection over some (non-isotropic) line, and a translation parallel to this line.
    We call a reflection over a non-isotropic line an \emph{axial symmetry}.
    The coset of $\SF$ contains, in particular, the set of axial symmetries.
    
    We define axial symmetries relative to non-isotropic lines only: if $\ell$ is non-isotropic, then $a$ being symmetric to $a'$ relative to $\ell$ means that $a-a'$ is normal to $\ell$, and also that for any $b\in \ell$, $d(a,b) = d(a',b)$.
    However, if $\ell$ is isotropic, this means that $a,a'$ lie on $\ell$.
    
    For $x,y\in \F^2$, we write $x\sim_\ell y$ to mean that $x$ is axially symmetric to $y$, relative to the (non-isotropic) line $\ell$. 
    
    The composition of two axial symmetries, relative to distinct lines $\ell$ and $\ell'$, is generally a rotation around the axes intersection point, by twice the angle between the lines. This is the same as in the Euclidean setting. If the lines are parallel, it is a translation in the normal direction (note that $\ell,\ell'$ are non-isotropic lines).
    
    In the sequel, for convenience of working within the group structure of $\SF$, rather than its other coset, we map the set of all axial symmetries into the group $\SF$.
    We map an axial symmetry to $\SF$ by composing it with the fixed axial symmetry $\rho$ relative to a non-isotropic line $\ell_\tau$. \label{elltau}
    
    The image of the set of axial symmetries under this mapping is the set of rotations around all points on $\ell_\tau$, which we denote by $R_\tau$.
    
    If $\ell_\tau$ is the $x$-axis, then explicitly
    \[
        R_\tau = \left\{\begin{pmatrix}
            u&-v& x_0(1-u)\\
            v&u&-x_0v\\
            0&0&1
        \end{pmatrix} \colon u^2 + v^2 = 1, u,v,x_0\in \F
        \right\}\,.
    \]
    A short calculation shows that, for this choice of $\ell_\tau$, the image of $R_\tau$ under the kinematic mapping is contained in the plane $X_2=0$.
    By Proposition~\ref{prop:kinematic-equivariant}, we see that $R_\tau$ is contained in a plane for any choice of $\ell_\tau$.
    This transformation motivates the role of incidence geometry.

\subsection{Incidence geometry}\label{sec:incidences}
  The key tool that we will use to estimate $\Dpin(A)$ is an incidence bound between points and planes in $\PF$ by the fourth author \cite{rudnev2014number}; for a selection of applications of this bound, see the survey \cite{rudnev2018point-plane}.
  
    \begin{thm}[Points-Planes in $\PF$] \label{t:rudnev}
        Let $\mathcal{P}$ be a set of points in $\F^3$ and let $\Pi$ be a set of planes in $\PF$, with $|\P|\leq |\Pi|$. 
        If $\mathbb{F}$ has positive characteristic $p$, suppose that  $|\P|\ll p^2$. 
        Let $k$ be the maximum number of collinear points in $\mathcal{P}$. Then
        \[
            \I(\mathcal{P},\Pi)\ll |\P|^{1/2}|\Pi| +k |\Pi|.
        \] 
    \end{thm}
  
  We also require a version of the incidence bound for restricted incidences. For a finite set of lines $\L\subseteq\F^3$ we define for any finite sets of points $P$ and planes $\Pi$ the number of ($\L$-)restricted incidences as
    
    \[
        \I_\L(P,\Pi):=|\{(p,\pi)\in P\times\Pi\colon \mbox{$p\in\pi$ and $\forall\ell\in \L$, $p\not\in \ell$ or $\ell\not\subseteq\pi$}\}|.
    \]

    \begin{thm}
    \label{thm:pt-plane-restricted}
    Let $\F$ be a field of characteristic $p$.
          Given a set of points $P\subseteq \F^3$ and sets of lines $\L$ and planes $\Pi$ contained in $\PF$, let $\mu$ denote the maximum number of points of $P$ (planes of $\Pi$) incident to (containing) a certain line not in $\L$.
    If $|P|=|\Pi|=N\ll p^2$, then
    \[
    \I_\L(P,\Pi)\ll N^{3/2}+\mu N.
    \]
    \end{thm}

    The proofs of the main results proceed by first relating the quantity $\Dpin(A)$ to the count of isosceles triangles. We count the number of isosceles triangles by studying the set of bisectors determined by $A$ considering the axial symmetries relative to this line set. Using the \BG{} embedding, we rephrase this as an incidence bound between points and planes. 

\subsection{Isosceles triangles and bisector energy}\label{sec:bisect-energy-isosc}

    The connection between $\Dpin(A)$ and the number of isosceles triangles with all vertices in $A \subset \R^2$ goes back to at least an argument of Erd\H os \cite{erdos1975some}. The same argument can be made to work over general fields. It is advantageous to exclude isosceles triangles where the two equal length sides have zero length. This decreases $\Dpin(A)$ by only one element (see Lemma~\ref{lem:pin-iso-lower-1} below) while it allows for much needed flexibility, for example when applying the kinematic mapping of \BG{}. In other words we are interested in bounding the number of isosceles triangles of the form
    \[
    \{ (a,b,c) \in A^3 \colon d(a,b) = d(a,c) \neq 0 \}.
    \]
    For such triangles the perpendicular bisector $\bisector(b,c)$ is non-isotropic and hence $b-c$, which is orthogonal to $\bisector(b,c)$ by the definition of $\bisector(b,c)$, is also non-isotropic
    . Hence the count of isosceles triangles of the above form is precisely
    \[
	\mathcal{T}_{\text{NI}}(A):=|\{ (a,b,c) \in A^3 \colon d(a,b) = d(a,c) \neq 0, b-c \text{ is non-isotropic} \}|.
    \]
    (The subscript NI is intended to remind the reader that $\mathcal{T}_{\text{NI}}$ is the count of non-isotropic
     isosceles triangles with vertices in $A$.)
    Our methods only require that the vector $b-c$ (equivalently $\bisector(b,c)$) is non-isotropic. We therefore define \emph{non-degenerate isosceles triangles} as isosceles triangles with non-isotropic base:
    \begin{align*} 
    \T(A) & = |\{ (a,b,c) \in A^3 \colon d(a,b) = d(a,c), b-c \text{ is non-isotropic} \}| \\
    	     & = |\{ (a,b,c) \in A^3 \colon d(a,b) = d(a,c), \bisector(b,c) \text{ is non-isotropic} \}|.
    \end{align*}
    
        The number of non-degenerate isosceles triangles determined by $A$ is inversely proportional to the number of pinned distances determined by $A$.
    \begin{lem}
      \label{lem:pin-iso-lower-1}
    If $A$ is a subset of $\F^2$ with at most $M$ points on a line, then
    \begin{equation}\label{eq:T-pind-cs}
	    |A|(|A|-2M+1)^2\leq (\Dpin(A)+1)(\T(A)+|A|^2)\,.   
    \end{equation}
    \end{lem}
    Lund and the second author prove this lemma as part of the proof of their Theorem~1.1~\cite[Section 2.5]{lund2018pinned}; we provide the proof here since it is fundamental.    \begin{proof}
    Let $C_r$ be the set of points in $\F^2$ of distance $r$ from the origin, and denote by $\Delta(A,a)$ the set of non-zero distances determined by $a$.
    Then
    \[
    |A|(|A|-2M+1)\leq 
    \sum_{a\in A}\sum_{r\in \Delta(A,a)}|A\cap (a+C_r)|,
    \]
    so by Cauchy-Schwarz and the bound $|\Delta(A,a)|\leq\Dpin(A)+1$, valid for any $a\in A$, we have
    \[
    |A|(|A|-2M+1)^2\leq (\Dpin(A)+1) \left( \sum_{a\in A}\sum_{r\not=0} |A\cap (a+C_r)|^2\right).
    \]
    We have
    \begin{equation}
      \label{eq:5}
      \sum_{a\in A}\sum_{r\not=0} |A\cap (a+C_r)|^2= \mathcal{T}_{\text{NI}}(A)+|A|^2,
    \end{equation}
    since the sum on the left hand side is equal to $\mathcal{T}_{\text{NI}}(A)$ plus the number of triples $(a,b,b')$ in $A^3$ such that $d(a,b)=d(a,b')\not=0$ and $b-b'$ is isotropic.
    If $b\not=b'$, then it's easy to conclude that $d(a,b)=d(a,b')=0$ (see e.g. \cite[Lemma 2.3]{lund2018pinned}). This is a contradiction, so there are $|A|^2$ such triples.    Now apply the bound $\mathcal{T}_{\text{NI}}(A) \leq \T(A)$.
    \end{proof}
    
    In order to bound $\T(A)$, we first introduce some terminology: let $(a,b,c)\in A^3$ be the vertices of a (non-degenerate) isosceles triangle determined by $A$. We call $a$ the \emph{apex} of the triangle, $b,c$ the \emph{base pairs}, and $\bisector(b,c)$ the \emph{symmetry axis} of the triangle. It may of course be the case that the isosceles triangle is in fact an equilateral triangle; in this situation it does not matter which of the three vertices is identified as the apex of the triangle. 
    
    If $(a,b,c)$ is a triangle contributing to the count of $\T(A)$ with apex $a$, then the perpendicular bisector of $b$ and $c$ passes through $a$. 
    Hence we can count the number of triangles by studying the perpendicular bisectors of pairs of points of $A$. 
    This notion is not novel, but we recall this principle in order to assist the reader to parse our notation.
    
    Given a line $\ell$, let $i_A(\ell)=|A\cap \ell|$ denote the number of points of $A$ incident to the line $\ell$, and let $b_A(\ell)$ denote the number of ordered pairs of distinct points in $A$ whose perpendicular bisector is $\ell$. We will use a modification of $b_A$ that suits the definition of $\T(A)$:
    \[
b_A^*(\ell) = |\{(b,c) \in A \times A  : b \neq c, \bisector(b,c) = \ell \text{ and $b-c$ is non-isotropic}\}|. 
    \]
    Because $\bisector(b,c)$ is isotropic precisely when $b-c$ is isotropic, we have $b_A^*(\ell) =0$ for all isotropic lines $\ell$, and that $b_A(\ell) = i_A(\ell)^2 - i_A(\ell)$ for all isotropic lines. Therefore 
    \begin{equation} \label{eq:S0}
    \sum_{\ell} b_A^*(\ell) = \sum_{\ell} b_A(\ell) - \sum_{\text{isotropic } \ell} \left( i_A(\ell)^2 - i_A(\ell)\right) = |A|^2 - |A| - |S_0|,
    \end{equation}    
    where $S_0$ is defined in \eqref{eq: Sr} on p.~\pageref{eq: Sr}.
    
    The number of non-degenerate isosceles triangles is then precisely the number of possible apexes of triangles and the number of non-isotropic base pairs:
    \begin{equation}
      \label{eq:T-sum}
      \T(A)=\sum_{\ell}i_A(\ell)b_A^*(\ell).
    \end{equation}

    The \emph{bisector energy}\/ of a set $A$ in $\F^2$ is the second moment of $b_A$. 
    \[
        \B(A):=|\{(a,b,a',b')\in A^4\colon \bisector(a,b)=\bisector(a',b')\}|=\sum_\ell b_A(\ell)^2.
    \]
    We write $\B^*(A)$ for the second moment of $b_A^*(\ell)$; this modified bisector energy allows us to avoid pathologies arising from isotropic vectors.
    
    Our next lemma bounds $\T(A)$ in terms of $\B^*(A)$. 
    \begin{lem}  \label{lem:bisectors-control-triangles}
        If $A$ is a subset of $\F^2$, then
        \[
            \T(A) \leq 2 |A|\B^*(A)^{1/2}.
        \]
    \end{lem}
    \begin{proof} 
       By \eqref{eq:T-sum} and an application of Cauchy-Schwarz, 
       \begin{align*}
       \T(A)  & \leq \sum_{\ell:i_A(\ell)=1}b_A^*(\ell) + \left( \sum_{\ell: i_A(\ell) >1} i_A(\ell)^2 \right)^{1/2} \left( \sum_\ell b_A^*(\ell)^2 \right)^{1/2} \\
       		& \leq \sum_{\ell}b_A^*(\ell) + |A| \B^*(A)^{1/2}.
       \end{align*}
         By \eqref{eq:S0}, we have
          \begin{align*}
          \sum_{\ell}b_A^*(\ell) & \leq \sum_{\ell}b_A^*(\ell)^{3/2} 
          				 \leq \left(\sum_{\ell}b_A^*(\ell)^2\right)^{1/2} \left(\sum_{\ell}b_A^*(\ell)\right)^{1/2} 
				   	 \leq |A| \B^*(A)^{1/2}\,.
	\end{align*}
    \end{proof}

\section{Proof of Theorem~\ref{thm: Pinned F}}\label{sec:proving pinned F}
In this section we prove the following bound on $\T(A)$, which is half of Theorem~\ref{thm: Pinned F}.
 \begin{thm}
      \label{thm:isosceles F}
        Let $\F$ be a field of characteristic $p$.
        If $p>0$ suppose $A\subseteq\F^2$ has cardinality $|A|\leq p^{4/3}$, then the number of non-degenerate isosceles triangles determined by $A$ satisfies
        \[
            \T(A) \ll |A|^{7/3}.
        \]
  \end{thm}

\subsection{Proof of Theorem~\ref{thm:isosceles F}}  \label{sec: 4.1}
   Our main technical result bounds the bisector energy of $A$ in terms of $|S_r(A)|$, defined in \eqref{eq: Sr} on p.~\pageref{eq: Sr}.
    \begin{prop}
      \label{prop:bisectors}
        Let $\F$ be a field of characteristic $p$.
        If $p>0$ suppose that $A\subseteq\F^2$ has cardinality $|A|\leq p^{4/3}$ and let $M$ denote the maximum number of collinear or co-circular points of $A$.
        Then the bisector energy of $A$ satisfies
        \[
            \B^*(A)\ll M|A|^2 + \sum_r |S_r|^{3/2}\,.
        \]
    \end{prop}

 \begin{proof}
    Recall that $S_r\subseteq A\times A$ denotes the set of segments of length $r$, and $C_r$ denotes the set of points in $\F^2$ of distance $r$ from the origin.
    We have 
    \[
        |S_r|=\sum_{a\in A}|A\cap (a+C_r)|.
    \]

    Let $\mathrm{Ax}_{(c,d)}$ be the set of elements $(x,y)\in \F^2\times \F^2$ that are axially symmetric to $(c,d)\in \F^2$ (with respect to some non-isotropic line).
    For a set $X\subseteq A\times A$ containing no isotropic segments (that is, $d(a-b, 0)\not=0$ for all $(a,b)\in X$), let $\mathcal{A}(X) := \{\mathrm{Ax}_x \colon x\in X\}$ consist of sets of elements attainable from $X$ via axial symmetries, partitioned according to the $x\in X$ characterising the axial symmetry. 
    Then, letting $\calL_\perp(A)$ denote the set of non-isotropic perpendicular bisectors of $A$, we have
    \begin{align*}
    \B^*(A) &=\sum_\ell b_A^*(\ell)^2\\
       &=|\{(a,b,c,d,\ell)\in A^4\times \calL_\perp(A) \colon (a,b)\sim_\ell (c,d)\}|\\
       & = \sum_{r\neq 0}|\{((a,b),(c,d))\in S_r^2 \colon (a,b)\in\mathrm{Ax}_{(c,d)}\}| \\
       & \quad \quad \quad + |\{((a,b),(c,d))\in S_0^2 \colon \exists~\ell, ~(a,b)\sim_\ell(c,d)\}|\\
       &=\sum_{r\neq 0}\I(S_r, \mathcal{A}(S_r)) + \calE,
     \end{align*}
    where the last line is a definition.
    
    The term $\calE$ can be bounded by $2M|A|^2$, since for each $a$, there are at most $2M$ choices of $b$ such that $a-b$ is isotropic; if $a,b$ and $c$ are chosen so that $c$ is axially symmetric to $a$, then there is only one choice for $d$, which gives the claimed bound.
    
    Let us prove the following claim, which immediately proves Proposition~\ref{prop:bisectors}.
    \begin{claim}
        \label{claim}
        Let $r\neq 0$, and suppose, if $\F$ has positive characteristic $p$, that $|A|\leq p^{4/3}$.
        Suppose that at most $M$ points of $A$ are collinear or co-circular in $\F^2$.
        Then
        \[
            \I(S_r,\mathcal{A}(S_r))\ll M|S_r| + |S_r|^{3/2}.
        \]
    \end{claim}

    \begin{proof}[Proof of Claim~\ref{claim}]
    The value of $\I(S_r,\mathcal{A}(S_r))$ is the same if we work in an extension of $\F$, and so we assume without loss of generality that $\F$ is algebraically closed. 
        
        Let us   embed the set $S_r$ in $\SF$ as follows. First we fix an element $(a,a') \in S_r$. Then we identify the segment $s_r\in S_r$  with the inverse of the rigid motion that takes $s_r$ to $(a,a')$.
        This rigid motion always exists, for one can translate $a$ to the origin, and then find the corresponding rotation, for $r\neq 0$. 
        Let $G_r$ denote the set of transformations in $\SF$ corresponding to segments in $S_r$. To summarise, the point $g \in G_r$ corresponds to a fixed line segment $s_r \in A^2 $ of length $r$.

        To each segment in $S_r$, we will associate a \emph{projective plane}\/ in \PF{}. More precisely, the plane is associated to the set $\mathrm{Ax}(s_r)$ characterised by a segment $s_r \in S_r$. 
        We proceed to describe this.  
        
        We will first choose, following the notation of p.~\pageref{elltau}, a non-isotropic line  $\ell_\tau$ so that for all $g,h\in G_r$, the transformation $g^{-1}h$ has no fixed points on  $\ell_\tau$; this is possible since $\F$ is algebraically closed, so there are infinitely many choices of $\ell_\tau$, but only a finite number of products $g^{-1}h$. By a slight abuse of notation, define also $\tau \in \SF$ to be the map of reflection over the fixed line $\ell_\tau$.
        Recall that $R_\tau$ is the set of axial symmetries composed with a reflection about $\ell_\tau$ and that $\kappa(R_\tau)$ is contained in a projective plane, which we also denote by $\kappa(R_\tau)$. We remind the reader that $\kappa$ is the kinematic mapping. 
        
        For $g \in \SF $ a transformation corresponding to an element in $S_r$, it follows from Proposition~\ref{prop:kinematic-equivariant} that the transformation $\phi_g$ is projective. 
        In particular, since $\kappa(R_\tau)$ is the entire plane (where $\kappa $ is the kinematic mapping) it follows from projectivity that the set $\phi_g(\kappa(R_\tau))$ is a projective plane in $\PF$. 
         
        Let $\Pi=\{\phi_g(\kappa(R_\tau))\colon g\in G_r\} = \{ \kappa(g R_\tau) \colon g \in G_r\}$.  
        We have $|\Pi|=|G_r|=|S_r|$. Indeed,  $\phi_g(\kappa(R_\tau))=\phi_h(\kappa(R_\tau))$ if and only if $g^{-1}h\in R_\tau$; on the other hand, every element of $R_\tau$ fixes a point on $\ell_\tau$, while no product $g^{-1}h$ with $g$ and $h$ in $G_r$ fixes a point on $\ell_\tau$. Hence there are no collisions of the form $\phi_g(\kappa(R_\tau))=\phi_h(\kappa(R_\tau))$. To summarise, we have mapped a line segment in $S_r$ to a plane in $\Pi$ so that different segments yield different planes. In this way, we map $\mathcal{A}(S_r)$ to a set of planes. 
        
        An identical process maps $S_r$ to points. First, we identify $s_r \in S_r$ with its axially symmetric counterpart; we map the symmetric segment to an element of $\SF$  and then we apply the kinematic mapping $\kappa$: 
        let $G_r'$ denote the set of transformations (obtained as in the definition of the set $G_r$) in $\SF$ corresponding to segments of length $r$ in $\tau(A) \times \tau(A)$ and set $P =  \kappa(G_r')$. Note that $\tau(A)$ is the reflection of $A$ in the line $\ell_\tau$. Since the kinematic mapping is injective, we have that $|P| = |G_r'| = |S_r|$.

       Note that the same procedure was applied to $S_r$ and $\mathcal{A}(S_r)$ to obtain $P$ and $\Pi$ respectively. 
       
      We will show that 
       \[
           \I(S_r,\mathcal{A}(S_r)) = \I(P,\Pi).
         \]

        Suppose that $\pi=\phi_g(\kappa(R_\tau))$ for some $g \in G_r$ and $\mathtt{q}=\kappa(h)$ for some $h\in G_r'$.  
        
        If we have an  incidence $\mathtt{q}\in\pi$, this is precisely the statement  $\kappa(h)\in  \phi_g(\kappa(R_\tau))$. Since $\SF$ is a group, it follows directly from Proposition~\ref{prop:kinematic-equivariant} that $\phi_g^{-1} = \phi_{g^{-1}}$. We thus obtain
        \[
            \kappa(g^{-1}h)=\phi_g^{-1}(\kappa(h))\in\kappa(R_\tau),
        \]
        so $h\in gR_\tau$.
        Now, let $(a,a')$ correspond to $g$ and $(b,b')$ correspond to $h$, so that
        \[
            g(a,a')=h(\tau(b),\tau(b'))=s_r.
        \]
        Since $h\in gR_\tau$, we have
        \[
            (\tau(b),\tau(b'))\in R_\tau^{-1}(a,a')=R_\tau(a,a'),
        \]
        thus $(b,b')$ is attainable from $(a,a')$ by an axial symmetry. Hence $           \I(S_r,\mathcal{A}(S_r)) = \I(P,\Pi)$ as claimed.

        We apply Theorem~\ref{t:rudnev} to $P$ and $\Pi$, claiming that the number of collinear  points is bounded by $M$. 

       This claim is is a direct consequence of \cite[Lemma 2.2]{lund2018pinned}.
       Firstly, note that points in $P:= \kappa(G_r')$ correspond to a segment $s_r \in S_r$ in the configuration space; $s_r$ determines the rigid motion $g$ via the relation $g^{-1}(\tau(s_r)) = (a,a')$ and $\kappa(g)$ is a point in $P$. This process is reversible.

        Lund and the second author \cite[Lemma 2.2]{lund2018pinned} show that if we have two segments $s,s' \in S_r$ , then any segment $s''\in S_r$ that is axially symmetric to both $s$ and $s'$ has endpoints lying on \emph{either} a pair of concentric circles \emph{or} a pair of parallel lines. The concentric circles, resp. parallel lines, are uniquely defined by $s$ and $s'$.  Moreover, the endpoints of both $s$ and $s'$ also lie on the concentric circles, resp. parallel lines. Via the chain of (bijective) mappings described above that enable us to obtain $P$ from $S_r$, one can calculate that segments lying on concentric circles (resp. parallel lines) are mapped to collinear points.  
        
        Thus
       \[
            \I(S_r,\mathcal A(S_r) ) = \I(P,\Pi) \ll M|S_r|+|S_r|^{3/2}\,.
        \]
        If $\F$ has positive characteristic $p$, then also we need the estimate $|\Pi|\ll p^2$; since $|\Pi|=|S_r|\ll |A|^{3/2}$ by Erd\H{o}s' bound on the number of times a distance can repeat \cite{erdos1946distances}, we have the required constraint for $|A|\ll p^{4/3}$.
        
       This completes the proof of the Claim. 
    \end{proof}
 Using $\sum_r|S_r|\leq |A|^2$, completes the proof of Proposition~\ref{prop:bisectors}.
\end{proof}

We now conclude the proof of Theorem~\ref{thm:isosceles F}. From Proposition~\ref{prop:bisectors} and Lemma \ref{lem:bisectors-control-triangles}, we have, with $M$ the maximum number of collinear or co-circular points of $A$, that 
\begin{equation}
    \label{eq:estimating T* with M}
    \T(A)\ll |A|\left(M|A|^2+\sum_{r\neq 0}|S_r|^{3/2}\right)^{1/2}\ll M^{1/2}|A|^2+|A|^{3/2}\left(\sum_{r\neq 0}|S_r|^2\right)^{1/4}\,.
\end{equation}

In the above, we use Cauchy-Schwarz as well as the bound $\sum_r|S_r|\leq |A|^2$.

Observe that $\sum_r|S_r|^2$ can be bounded in terms of $\T$. Indeed,
\begin{align*}
    \sum_{r\neq 0}|S_r|^2 & = \sum_{r\neq 0}\left(\sum_{a\in A}\sum_{b\in A}\one_{\|a-b\| = r}\right)^2\\
    & \leq|A| \sum_{r\neq 0} \sum_{a\in A}\sum_{b,b'\in A}\one_{\|a-b\| =r}\one_{\|a-b'\| = r}\le |A| (\T(A) + 3|A|^2)\,,
\end{align*}

where we used the fact that the number of isosceles triangles $(a, b, b')$ with $||b-b'||=0$ is at most $3|A|^2$ because for any $a,b \in A$, there are at most 3 admissible $b'$ lying in the intersection of a circle of non-zero radious centered at $a$ and one of two isotropic lines through $b$.

Either $\T(A) \le |A|^2$, in which case the theorem is proved, or $\T(A) + 3|A|^2 \ll \T(A)$. Hence we may assume
\[
    \T(A)\ll M^{1/2} |A|^2 + |A|^{7/4}\T(A)^{1/4}\,.
\]

If the second term dominates, then $\T(A)\ll |A|^{7/3}$, as required. 
If the first term dominates, then $\T(A)\ll \min(M^2|A|,M^{1/2}|A|^2)$. When $M\leq |A|^{2/3}$ then the first minimand ensures that $\T(A)\ll |A|^{7/3}$. When $M\geq |A|^{2/3}$, we use an additional argument to obtain the statement of Theorem~\ref{thm:Isosc Triang}, which is independent of $M$. We remark that $\Dpin(A)\gg M$ (when $A$ is not contained in an isotropic line), so this additional argument is not necessary to prove Theorem~\ref{thm: Pinned F}. 

\subsection{Additional argument to remove rich circles and lines}
\begin{lem}[Pruning heavy circles and lines I]
      \label{lem:3}
        Suppose that $A$ is the disjoint union of $B$ and $C$.
        If all of the points of $C$ are contained in a circle or a line, then
        \[
            \T(A)\leq \T(B) + 8|A|^2.
        \]
\end{lem}
\begin{proof}
$\T(A) - \T(B)$ is the number of non-degenerate isosceles triangles with at least one vertex in $C$. Denoting by $(x,y,z)$ an isosceles triangle with apex $x$ and base pairs $(y,z)$ we see that every non-degenerate triangle contributing to $\T(A) - \T(B)$ must be either of the form $(c,a,a')$, $(a,c,a')$ or $(a,a',c)$ for some $c \in C$ and $a,a' \in A$. In the first case $c$ belongs to $C \cap \bisector(a,a')$ and in both the second and third cases $c$ belongs to $C \cap \sigma_a(d(a,a'))$, where $\sigma_x(\rho)$ is the circle of radius $\rho$ centred at $x$.

Let us first prove the lemma when $C$ is contained in a line. In the first case $|C \cap \bisector(a,a')| \leq 1$ unless $C \subseteq \bisector(a,a')$, in which case it is at most $|A|$. Hence for the at most $|A|^2$ pairs $(a,a') \in A^2$ for which $C \not\subseteq \bisector(a,a')$ there is at most one $c$ so that $(c,a,a')$ contributes to $\T(A) - \T(B)$. Since $\bisector(a,a')$ is non-isotropic there are at most $|A|$ pairs $(a,a') \in A^2$ such that $\bisector(a,a') \supseteq C$. Therefore the total contribution to $\T(A) - \T(B)$ from such base pairs is at most $|A|^2$.

To treat the second and third cases, we note that $|C \cap \sigma_a(d(a,a'))| \leq 2$ for all $(a,a') \in A^2$ with $d(a,a') \neq 0$. Hence the total contribution to $\T(A) - \T(B)$ from such base pairs is at most $2|A|^2$. 
We must also consider the case when $d(a,a')=0$. For this situation, we fall into two cases: 

{\bf Case $1$:} Assume $C$ is supported on a non-isotropic line $\ell_C$. We know that  $\sigma_a(d(a,a'))$ is an union of two isotropic lines, say $\ell_1$ and $\ell_2$, assume $\ell_{1}$ contains $a$ and $a'$. The triangles $(a,c,a')$ or $(a,a',c)$ can contribute to $\T(A) - \T(B)$ only if  $c\not\in \ell_{1}$, otherwise we have $||c-a||=||c-a'||=0$, contradicting the non-degeneracy condition. So, $c$ is determined uniquely by the intersection of $\ell_2$ and $\ell_C$. 

{\bf Case $2$:} Assume $C$ is supported on an isotropic line $\ell_C$ of slope $i$ (the case of slope $-i$ is treated in the same way). A triangle $(a,c,a')$ or $(a,a',c)$ contributes to $\T(A)$ only if $\|a'-c\| \neq 0$. 
If $a \notin \ell_C$ and $\|a-a'\| = 0$, then there is at most one triangle $(a, c, a')$ or $(a, a', c)$ with $||a-c||=0$. 
If $a\in \ell_C$ then 
it is not hard to see that the number of triangles with $a\in \ell_C$ is at most $|\ell_C\cap A|\cdot |A|\le |A|^2$. 


In other words, we have $\T(A) - \T(B) \leq 8|A|^2$. 

Suppose now that $C$ is contained in a circle of non-zero radius. We proceed in a very similar fashion and omit some details. The number of $(c,a,a')$ that contribute to $\T(A) - \T(B)$ is at most $2 |A|^2$ (at most two $c$ work for each $(a,a')$). The number of $(a,c,a')$ with $C \not\subseteq \sigma_a(d(a,a'))$ is at most $2|A|^2$ (at most two $c$ work for each $(a,a')$). The number of $(a,a')$ with $C \subseteq \sigma_a(d(a,a'))$ that contribute to $\T(A) - \T(B)$ is at most $|A|$ ($a$ is fixed and $a'$ is free), and for each such pair there are at most $|A|$-many c that are admissible. In total $\T(A) - \T(B) \leq 8|A|^2$. 
    \end{proof}
    
    \begin{lem}[Pruning heavy lines and circles II]
      \label{lem:1}
        There is a subset $A'\subseteq A$ such that at most $|A|^{2/3}$ points of $A'$ are collinear or co-circular, and
        \[
            \T(A)\leq \T(A') + 8 |A|^{7/3}.
        \]
    \end{lem}
    \begin{proof}
        Use Lemma~\ref{lem:3} to greedily remove lines and circles, gaining a factor of $8 |A|^2$ each time.
        If we only remove lines and circles with more than $|A|^{2/3}$ points on them, then this procedure terminates after $|A|^{1/3}$ steps.
    \end{proof}
    
    \subsection{Conclusion of the proofs of Theorem~\ref{thm: Pinned F} and \ref{thm:isosceles F}}
    
    For Theorem~\ref{thm:isosceles F} we apply Lemma~\ref{lem:1} to get the $A' \subseteq A$ for which the number $M'$ of colinear or cocircular points is at most $|A|^{2/3}$ and   
    \[
         \T(A) \ll \T(A') + |A|^{7/3}.
    \]
By the work in Section~\ref{sec: 4.1}, summarised in its concluding paragraphs, we have
\[
\T(A') \ll |A'|^{7/3} + (M')^{1/2} |A'|^2 \ll |A|^{7/3},
\]
which completes the proof of Theorem~\ref{thm:isosceles F}.

For Theorem~\ref{thm: Pinned F} we first note that because $A$ is not contained in an isotropic line we have $\Dpin(A) \geq (M-1)/2$. Therefore, if $M > |A|/3$, then $\Dpin(A) \gg |A|$. Otherwise, we apply Lemma~\ref{lem:pin-iso-lower-1} and simply note that $|A| - 2M+1 \gg |A|$:
\[
\Dpin(A) \gg \frac{|A|^3}{\T(A)} \gg |A|^{2/3}.
\]

\section{Conclusion of the proof of Theorem~\ref{thm:Isosc Triang} and proofs of Theorems~\ref{thm:Pinned Falconer}~and~\ref{thm: Pinned Fp only}}

As before, the proof of the pinned distance statements of Theorems~\ref{thm:Pinned Falconer}~and~\ref{thm: Pinned Fp only} rely on a new estimate for the number of non-degenerate isosceles triangles with vertices in $A$. Our final important task is to prove the following.

   \begin{thm}\label{thm:Isosc Triang F_p}
        Let $p$ be a prime and $A \subseteq \F_p^2$. If $p \leq |A| \leq p^{4/3}$, then
        \begin{equation}\label{eq:T*bound F_p}
            \T(A) = \frac{|A|^3}{p} + O(p^{2/3} |A|^{5/3} + p^{1/4}|A|^2).
        \end{equation}
    \end{thm}

  We prove Theorem \ref{thm:Isosc Triang F_p} using a modification of the proof of Theorem~\ref{thm:isosceles F}. 
    The improvement does \emph{not} come from an improved bound on the bisector energy $\B^*(A)$. Indeed, Lund, Sheffer and de Zeeuw  \cite[Theorem 2.2]{lund2016bis} gave an example of a point configuration with $\B^*(A)$ being forbiddingly large, owing to the $M$ quantity above. Our strategy is of a structural nature which lends itself to a dichotomy argument: we show that the components of sets that lead to large $\B^*(A)$ cannot supply many isosceles triangles. Before giving the proof let us quickly derive Theorems~\ref{thm:Isosc Triang}~and~\ref{thm:Pinned Falconer}~and~\ref{thm: Pinned Fp only}.

The first statement in Theorem~\ref{thm:Isosc Triang} follows by combining Theorem~\ref{thm:Isosc Triang F_p} (for the range $p \leq |A| \leq p^{4/3}$) with Theorem~\ref{thm:isosceles F} (for the range $ |A| < p$). The second statement in Theorem~\ref{thm:Isosc Triang} is Theorem~\ref{thm:isosceles F}. 
  
  The first statement in Theorem~\ref{thm:Pinned Falconer}~and~Theorem~\ref{thm: Pinned Fp only} follow from Theorem~\ref{thm:Isosc Triang F_p} and Lemma~\ref{lem:pin-iso-lower-1} by noting that $M \leq p\pm 1 \leq |A|/4$. (When $p\equiv 1 \mod(4)$, $M\leq p-1$, and when $p\equiv 3\mod(4)$, $M\leq p+1$.)
  
  We are left with proving the second statement in Theorem~\ref{thm:Pinned Falconer}. For $a \in A$, recall the definition of $\Delta(A;a)=|\Delta(A,a)|$ in \eqref{eq:Dax} on p.~\pageref{eq:Dax}; set
  \[
  \T_a(A) = |\{ (b,c) \in A^2 : d(a,b) = d(a,c) , b-c \text{ is non-isotropic} \}|
  \]
  so that 
  \[
  \T(A) = \sum_{a \in A} \T_a(A);
  \]
  similarly define $\mathcal{T}_{NI,a}(A)$, and set $z_a = | A \cap \sigma_a(0)|$ be the number of elements of $A$ lying at distance zero from $a$. Note that $z_a \leq 2p$ for all $a \in A$.
  
  The Cauchy-Schwarz inequality gives that, for all $a \in A$,
  
 \[
  \frac{(|A| - z_a)^2}{\Delta(A;a)}-(|A|-z_a-1) \leq \mathcal{T}_{NI,a}\,.
   \]
To see this, we have 
\[|A|-z_a=\sum_{r\in \Delta(A, a), r\ne 0} |A\cap (a+C_r)|,\]
where $C_r$ is the circle centred at the origin of radius $r$. Therefore, 
\begin{align*}
(|A|-z_a)^2&\le \Delta(A; a)\cdot \left(\sum_{r\in \Delta(A, a), r\ne 0}|A\cap (a+C_r)|^2\right)\\
&=\Delta(A; a)\cdot \bigg(\mathcal{T}_{NI,a}+(|A|-z_a-1)\bigg),
\end{align*}
where the term $(|A|-z_a-1)$ comes from isosceles triangles $(a, b, b')$ with $b=b'$. As in the proof of Lemma \ref{lem:pin-iso-lower-1}, we also note that there is no isosceles triangle $(a, b, b')$ with $||a-b||=||a-b'||\ne 0$, $b\ne b'$, and $||b-b'||= 0$.

  From this it follows that, with the ensuing constant 4 chosen to simplify the exposition, 
  \[
  \frac{(|A| - 4p)^2}{p}-|A| \leq \mathcal{T}^{*}_a(A).
  \]
  By making $p$ large enough we ensure that $\mathcal{T}^{*}_a(A) > (1 - \varepsilon^2 /2) |A|^2/p$ for all $a \in A$; and also that $\T(A) \leq (1 + \varepsilon^2/2) |A|^3 /p$.
  
  Now let $\delta$ be the proportion of elements of $A$ for which $\mathcal{T}^{*}_a(A) > (1 + \varepsilon) |A|^2/p$. Summing over $a \in A$ and using the established lower bound on $\mathcal{T}^{*}_a(A)$, we get
  \[
  (1 - \frac{\varepsilon^2}{2}) (1 - \delta) + \delta (1 + \varepsilon) <\frac{p}{|A|^3}\sum_{a\in A}\T_a(A)\leq 1 + \frac{\varepsilon^2}{2}.
  \]
  Some algebra gives $\delta \leq \varepsilon$. Therefore for at least $(1 - \varepsilon) |A|$ many $a \in A$ we have 
  \[
  \Delta(A;a) \geq \frac{(|A| - 4p)^2}{\T_a(A)} \geq \frac{(1-\varepsilon) |A|^2}{(1+\varepsilon) |A|^2/p} \geq (1-\varepsilon)^2 p \geq (1-2 \varepsilon) p.
 \]
 Changing $\varepsilon$ to $\varepsilon/2$ finishes the proof.
    
    \subsection{Preliminary lemma:  bounding rich lines and circles}
    In preparation of proving Theorem~\ref{thm:Isosc Triang F_p} we prove a standard elementary combinatorial lemma about rich lines or circles. A $k$-\emph{rich} line or circle with respect to a set $A$ is a line or circle containing at least $k$ elements of $A$.

    \begin{lem} \label{lem:Rich}
    Let $A \subseteq \F_p^2$ and $k \geq \sqrt{8 |A|}$. Let $\gamma_1, \dots, \gamma_n$ be the complete list of $k$-rich circles and lines with respect to $A$. Then $n \leq 2|A|/k$ and for each $i=1,\dots, n$ we have
    \[
    \left| (A \cap \gamma_i) \setminus \bigcup_{j \neq i} \gamma_j \right| \geq \frac{|A \cap \gamma_i|}{2}.
    \]
    \end{lem}
    
    \begin{proof}
    Consider a maximal collection of circles and lines $\gamma_1^*, \dots, \gamma_m^*$ such that for each $i=1,\dots,m$
    \[
    \left| (A \cap \gamma_i^*) \setminus \bigcup_{j =1}^{i-1} \gamma_j^* \right| \geq \frac{k}{2}.
    \]
    It follows that
    \[
    \frac{mk}{2} \leq \left| A \cap \bigcup_{i =1}^{n} \gamma_i^* \right| \leq |A|,
    \]
    and hence $m \leq 2 |A| /k$. 
    
    Next, consider a circle or a line not in the collection. It is incident to strictly fewer than $k/2$ elements of $A \setminus \bigcup \gamma_j^*$ (else it must be added to the collection) and incident to at most 2 elements of each $\gamma_j^*$. Therefore it is incident to strictly fewer than $k/2 + 2m \leq k/2 + 4 |A| / k \leq k/2 + k/2 =k$ points of $A$. We deduce that every $\gamma_i$ is included in this collection and so $n \leq m \leq 2|A|/k$. 
    
    Finally, note that
    \[
    \left| \gamma_i \cap \bigcup_{j \neq i} \gamma_j \right| \leq 2 (n-1) \leq \frac{4 |A|}{k} \leq \frac{k}{2} \leq \frac{|A \cap \gamma_i|}{2}.
    \]  
    \end{proof}

    \subsection{Counting isosceles triangles}
    We first remove isotropic segments from consideration and obtain a balanced count of non-degenerate isosceles triangles.
    From \eqref{eq:S0} on p.~\pageref{eq:S0} we have the equality
    \[
        \T(A) - \frac{|A|^3}{p} = \sum_\ell \left( i_A(\ell) - \frac{|A|}{p} \right) b_A^*(\ell) + \frac{|A|}{p} |S_0|  + \frac{|A|^2}{p}.
    \]
    
    We use the bound $|S_0| \leq 2 p |A|$ (for each element in $A$ there are at most $2p$ elements of $A$ on the isotropic lines incident to it) to get
    \begin{equation}
    \label{eq:balancedT*}
        \T(A) - \frac{|A|^3}{p} \leq 3 |A|^2 +  \sum_\ell \left( i_A(\ell) - \frac{|A|}{p} \right) b_A^*(\ell).
    \end{equation}
    The $3|A|^2$ term is smaller than $p^{1/4} |A|^2$, so we could ignore it from now on. 

    Let us now present an overview of the argument. We partition the above balanced count of isosceles triangles \eqref{eq:balancedT*} into three parts. This is achieved by partitioning the first set of lines with $b_A^*(\ell)>0$ into two disjoint sets $\L_1$ and $\L_2$ (the partition is described below). We then further decompose the set of isosceles triangles with symmetry axis in $\L_2$ into two parts, and estimate the size of each separately: those with the base pair $(b,c)$ on a circle or line that contains `many' elements of $A$, and those with the base pair $(b,c)$ on a circle or line that contains `few' elements of $A$.
    
\subsection{Decomposition of $\T(A)$}    
\label{sec:decomp-of-T*A}
To define $\L_1$ we apply Lemma~\ref{lem:Rich} with $k= \sqrt{8|A|}$. Let $\Gamma=\{\gamma_1,\ldots,\gamma_n\}$ be the complete set of $k$-rich circles and lines with respect to $A$, and for $\gamma\in\Gamma$, let $A_\gamma =A \cap \gamma$. 

Let $C$ denote the set of centres of circles in $\Gamma$ and $V$ the set of directions of lines in $\Gamma$. For $c \in C$, let $\Gamma_c$ denote the set of circles in $\Gamma$ centred at $c$ and let $A_c=\bigcup_{\gamma\in\Gamma_c}A_\gamma$. Similarly for $v \in V$, let $\Gamma_v$ denote the set of lines in $\Gamma$ with direction $v$ and let $A_v=\bigcup_{\gamma\in\Gamma_v}A_\gamma$.

Now let 
\begin{equation*} \label{eq:K}
	K = \frac{|A|^{4/3}}{p^{2/3}}
\end{equation*}
(it may be more helpful to view $K$ as a parameter to be chosen later) and define $C_1 \subseteq C$ to be the set of centres $c$ for which $|A_c| > K$ and $V_1 \subseteq V$ to be the set of directions $v$ for which $|A_v| > K$.

Define $\L_1$ to be the collection of lines with $b_A^*(\ell) >0$ that pass through some centre in $C_1$ or that are orthogonal to some direction in $V_1$; and define $\L_2$ to be the remaining lines with $b_A^*(\ell)>0$. 


Now for $i=1,2$ set $\T_{i,\text{bal}}(A)$ to be the balanced sum
\[
	\T_{i,\text{bal}}(A) = \sum_{\ell \in \L_i}  \left( i_A(\ell) - \frac{|A|}{p} \right) b_A^*(\ell).
\]
In this notation we have
\begin{equation} \label{eq:Ti}
	\T(A) - \frac{|A|^3}{p} \leq 3|A|^2 + \T_{1,\text{bal}}(A) + \T_{2,\text{bal}}(A).
\end{equation}

We bound $\T_{1,\text{bal}}(A)$ by a direct counting argument. For $\T_{2,\text{bal}}(A)$ we tweak the method developed in the proof of Theorem~\ref{thm:Isosc Triang} to be able to deal efficiently with sets with large $M$.

\subsection{Bounding $\T_{1,\text{bal}}(A)$: triangles interacting with rich lines and circles}

    An isosceles triangle $(a,b,c)\in A^3$ with apex $a$ contributing to the count of $\T_{1,\text{bal}}$ satisfies (i) $d(a,b) = d(a,c)\neq 0$; (ii) $b-c$ non-isotropic; (iii) the perpendicular bisector of $b$ and $c$ either passes through a centre in $C_1$ or is orthogonal to a direction in $V_1$.
      
    \begin{claim}
            \label{T1}
    \[
    	\T_{1,\text{bal}}(A) \ll \frac{|A|^3}{K}.
    \]
    \end{claim}
    
    \begin{proof}[Proof of Claim~\ref{T1}]
    We begin by bounding the number of centres in $C_1$ and the number of directions in $V_1$. From the second claim in Lemma~\ref{lem:Rich} we get
    \[
    \frac{|C_1| K}{2} \leq \frac{\sum_{c \in C_1} |A_c|}{2} \leq \left| \bigcup_{c \in C_1} A_c \right| \leq |A|.
    \] 
    Therefore $|C_1| \leq 2 |A| /K$ and similarly $|V_1| \leq 2 |A| / K$. 
    
    For each $c \in C_1$ set $\L_{1,c}$ be the set of lines in $\L_1$ that are incident to $c$; similarly define $\L_{1,v}$ for each $v \in V_1$ to be the set of lines in $\L_1$ with direction $v$. Note that 
    \[
    \L_1 \subseteq \bigcup_{c \in C_1} \L_{1,c} \cup \bigcup_{v \in V_1} \L_{1,v}
    \]
    (we do not claim that $\L_{1,c}$ and $\L_{1,v}$ are disjoint). Note that $\sum_{\ell\in\L_{1,c}} i_A(\ell)\leq |A|+p$ (each point of $A$ is incident to at most one line in $\L_{1,c}$, unless $c\in A$) and $\sum_{\ell\in\L_{1,v}} i_A(\ell)\leq |A|$ (each point of $A$ is incident to at most one line in $\L_{1,v}$). Since $b^*_A(\ell)\leq |A|$ we have
    \begin{align*}
    \T_{1,\text{bal}}(A) 
    &\leq \sum_{c\in C_1}\sum_{\ell\in \L_{1,_c}} i_A(\ell) b_A^*(\ell) + \sum_{v\in V_1}\sum_{\ell\in \L_{1,v}} i_A(\ell) b_A^*(\ell) \\
    &\leq (|C_1| + |V_1|) (|A|+p) |A| \\
    & \ll \frac{|A|^3}{K} + \frac{p|A|^2}{K}\\
    & \ll \frac{|A|^3}{K},
    \end{align*}
    because $|A| \geq p$. This completes the proof of Claim~\ref{T1}.
    \end{proof}

\subsection{Triangles not interacting with rich lines and circles}
We now turn to bounding $\T_{2,\text{bal}}(A)$. It is helpful to assume that $\T(A) - |A|^3 p^{-1} \geq 6 |A|^2$ (otherwise we are done because $|A|^2 \leq p^{2/3} |A|^{5/3}$ for all $A \subseteq \F_p^2$). We may also assume that $\T_{1,\text{bal}}(A) \leq \T_{2,\text{bal}}(A)$ (otherwise we are done by Claim~\ref{T1} by recalling $K = |A|^{4/3} p^{-2/3}$) and consequently $\T(A) - |A|^3 p^{-1} \leq 4 \T_{2,\text{bal}}(A)$. 

We will prove the following bound:
\begin{claim}
        \label{T2}
    If $\T(A) - |A|^3 p^{-1} \leq 4 \T_{2,\text{bal}}(A)$, then
    \[
    \T_{2,\text{bal}}(A) = O\left(p^{2/3} |A|^{5/3} + p^{1/4} |A|^2 + p^{1/2} |A|^{7/4} + p K^{1/2} |A|\right).
    \]
\end{claim}

We defer the (rather involved) proof of Claim~\ref{T2} and complete the proof of Theorem~\ref{thm:Isosc Triang}.

\subsection{Conclusion of the proof of Theorem \ref{thm:Isosc Triang} assuming Claim~\ref{T2}}
    We complete the proof of Theorem~\ref{thm:Isosc Triang F_p}. The value of $K$ that balances the first term in Claim~\ref{T1} and the last term in Claim~\ref{T2} is $K = |A|^{4/3} / p^{2/3}$, the value that was assigned to $K$ on page~\pageref{eq:K}. Combining \eqref{eq:Ti} with Claims \ref{T1} and \ref{T2} with $K = |A|^{4/3} / p^{2/3}$, and the paragraph above Claim~\ref{T2} yields
\[
\T(A) - \frac{|A|^3}{p} \ll p^{2/3} |A|^{5/3}  + p^{1/4} 	|A|^2 + p^{1/2} |A|^{7/4}.
\]
As $|A| \leq p^2$, the fourth summand is always smaller than the first. Since $|A| \geq p$, the second summand is also smaller than the first.

\section{Proof of Claim~\ref{T2}}
To prove Claim~\ref{T2} we follow a similar strategy to the proof of Theorem~\ref{thm:isosceles F}. We relate the contribution to $\T(A)$ to the \emph{bisector energy} of $A$; we bound this in the language of axial symmetries and finally, using the technology of the \BG{} mapping, we interpret this as an incidence bound between points and planes. This proof is better suited to the `large' $|A|$ case than Theorem~\ref{thm:isosceles F}.

    In Theorem~\ref{thm:isosceles F} we applied the standard formulation of the point-plane incidence theorem, resulting in a bound in terms of the number of collinear and cocircular points of $|A|$. We improve this using the version of the point-plane theorem for restricted incidences, Theorem~\ref{thm:pt-plane-restricted}; this motivates the definition of $\L_1$.

    We bound $\T_{2,\text{bal}}(A)$ in terms of the bisector energy using the Cauchy-Schwarz inequality:
    \begin{align*}
         \T_{2,\text{bal}}(A) 
        & =  \sum_{\ell \in \L_2} \left( i_A(\ell) - \frac{|A|}{p} \right) b_A^*(\ell) \\
        &\leq \left( \sum_\ell \left( i_A(\ell) - \frac{|A|}{p} \right)^2 \right)^{1/2} \left( \sum_{\ell\in \L_2} b_A^*(\ell)^2 \right)^{1/2}.
    \end{align*}
    The sum inside first bracket is at most $p|A|$ (see, for example, \cite[Lemma 1]{murphy2016point-line}). We denote the sum in the second bracket by $\B_2^*(A)$:
    \[
    \B_2^*(A) := \sum_{\ell \in \L_2}b_A^*(\ell)^2.
    \]
    In other words,
    \begin{equation} \label{T2*}
    \T_{2,\text{bal}}(A)  \leq \sqrt{p |A| \B_2^*(A)}.
    \end{equation}
    
    We are left with bounding $\B_2^*(A)$. It is worth noting that bounding the sum inside the first bracket by $p|A|$ rather than $|A|^2$ is one of the reasons why Theorem~\ref{thm:Pinned Falconer} applied to a set $A$ with $|A|= p^{5/4}$ yields $\Dpin(A) \gg |A|^{4/5}$ (better than the lower bound that Theorem~\ref{thm: Pinned F} gives).

    \subsection{Bounding bisector energy using axial symmetries}
    
    We first develop notation corresponding to the particular situation of $\T_{2,\text{bal}}(A)$.
    Let $\ax'_{(c,d)}$ be the set of elements $(x,y)\in\F_p^2  \times \F_p^2$ that are axially symmetric to $(c,d)$ with respect to some line in $\L_2$.
    Similarly, let $\mathcal{A}'(X):=\{\ax'_{(c, d)}\colon (c, d)\in X\}$.
    
    Using this language, we formulate a bound on $\B_2^*(A)$ as an incidence problem.
    
    \begin{align*}
        \B_2^*(A) 
		   & = \sum_{\ell\in \L_2}b_A^*(\ell)^2 \\
    	   &\le |\{(a,b,c,d,\ell)\in A^4\times\L_2 \colon\exists~\ell, ~ (a,b)\sim_\ell (c,d)\}|\\
      	   & \leq \sum_{r\neq 0}|\{((a,b),(c,d))\in S_r^2 \colon (a,b)\in\ax'_{(c,d)}\}| \\
	   & \quad  \quad \quad + |\{((a,b),(c,d))\in S_0^2 \colon \exists~\ell ,~ (a,b)\sim_\ell(c,d)\}|\\
      	   &=\sum_{r\neq 0} \I(S_r, \mathcal{A}'(S_r)) + \mathcal{E},
     \end{align*}
     where the last line is a definition. 
     
     To bound the (naturally insignificant) term $\calE$ we follow Hanson, Lund, and Roche-Newton~\cite[Proof of Theorem 3 from Lemma 10]{hanson2016distinct}. We show that
     \[
     \mathcal{E} \leq \sum_{r\neq 0}\I(S_r, \mathcal{A}'(S_r)) + 2 |A|^2
     \]
     and consequently that
     \begin{equation} \label{eq: E}
      \B_2^*(A) \ll \sum_{r\neq 0} \I(S_r, \mathcal{A}'(S_r)) + |A|^2.
     \end{equation}
     Suppose, without loss of generality that $\ell = \bisector(a,c) =  \bisector(b,d)$. It follows, by symmetry or other considerations, that $d(a,d) = d(b,c)$. If $d(a,d) = d(b,c) \neq 0$, then a permutation of the quadruple $(a,b,c,d)$ appears in some $\I(S_r, \mathcal{A}'(S_r))$ with $r \neq 0$. While if $d(a,d) = d(b,c) = d(a,b) = d(c,d) =0$ for every $(a,c) \in A^2$ both of $b,d$ must lie on isotropic lines through $a$ and through $c$. Let us call these lines $\ell_{a1}, \ell_{a2}, \ell_{c1}, \ell_{c2}$ The points $a-c$ and $b-d$ are non-isotropic and therefore these four lines must be distinct (otherwise either both $a$ and $c$ or $b$ and $d$ would be on an isotropic line). Therefore there are at most two possibilities for $(b,d)$: the two permutations of $\{\ell_{a1} \cap \ell_{a2}, \ell_{c1} \cap \ell_{c2}\}$.

\subsection{Application of restricted incidence bound}

       As in Theorem~\ref{thm:isosceles F} we work over $\F$, the algebraic closure of $\F_p$; $S_r$ and $\mathcal{A}'(S_r)$ embed into $\mathbb{FP}^3$ as sets of points $P$ and planes $\Pi$, respectively. The elements of $\mathcal{A}'(S_r)$ will be embedded into subsets of planes, which is sufficient for our purposes. We identify the subset of the plane that an element of $\mathcal{A}'(S_r)$ embeds into with the whole plane for ease of notation. 
       
As in the proof of Theorem~\ref{thm:isosceles F}, we will study incidences between points and planes. We map segments in $S_r$ to points $P$ in $\mathbb{FP}^3$ and we map pairs of segments in $S_r$ to planes $\Pi$ in $\mathbb{FP}^3$.

Suppose we have an incidence between a point $\mathtt{q}_x$ and distinct planes $\pi_y $, $\pi_z$ in $\mathbb{FP}^3$, where $\mathtt{q}_x$ corresponds to the segment $x = (x_1,x_2)\in S_r$ and the planes $\pi_y$ and $\pi_z$ correspond to the segments $y = (y_1,y_2)$ and $z=(z_1,z_2)$ in $S_r$  for $y\neq z$.

The incidence $\mathtt{q}_x \in \pi_y\cap \pi_z$ has the following interpretation in the original ``configuration space'' pertaining to the original set $A\subseteq \F_p^2$: the endpoints of the segment $x$ lie on two concentric circles or parallel lines $\gamma,\gamma'$; further, the same is true for the endpoints of $y$ and $z$: $x_1,y_1,z_1\in\gamma, x_2,y_2,z_2\in\gamma'$. The concentric circles or parallel lines are uniquely determined by $y$ and $z$.
This interpretation is described in Theorem~\ref{thm:isosceles F} and is a consequence of a lemma of Lund and the second author \cite[Lemma 2.2]{lund2018pinned}.

 We say that the pair of concentric circles $\gamma,\gamma'$ are the \emph{annulus} belonging to $\ax_y,\ax_z$. Every pair of planes in $\mathbb{FP}^3$ intersects in a line that defines an annulus as described in the previous paragraph. 

    Let $\L$ be the set of lines in $\mathbb{FP}^3$ determined by the intersection of two distinct planes in $\Pi$ such that these planes determine an annulus $(\gamma,\gamma')$ with $\gamma,\gamma'\in\Gamma$. The set $\Gamma$ is the complete set of $k$-rich circles and lines with respect to $A$ defined in Section~\ref{sec:decomp-of-T*A} on p. \pageref{sec:decomp-of-T*A}. Recall that $k=\sqrt{8|A|}$.
    
    For any pair of planes in $\Pi$ whose intersection $\ell$ is \emph{not} in $\L$, we have $|P\cap \ell| \ll |A|^{1/2}$. 

    We will apply the restricted incidence bound Theorem~\ref{thm:pt-plane-restricted} to this set-up. We may apply the theorem because, by Erd\H{o}s' argument \cite{erdos1946distances}, we have the bound $|S_r| \ll |A|^{3/2} \leq p^2$ (this is the only point we use the hypothesis $|A| \leq p^{4/3}$). Since $\mu\ll |A|^{1/2}$ and $N=|S_r|$ we get 
     \begin{equation} \label{eq: IL}
       \I_{\L}(S_r,\mathcal{A}'(S_r))\ll |S_r|^{3/2} + |A|^{1/2} |S_r|.
    \end{equation}

\subsection{Bounding $\I(S_r,\mathcal{A}'(S_r))$}
     For $r\not=0$, we decompose $\I(S_r,\mathcal{A}'(S_r))$ into two pieces.

    The quantity $\I(S_r,\mathcal{A}'(S_r))$ is bounded by the sum of $\I_\L(P,\Pi)$ and the number of pairs $(a,b),(a',b')\in S_r$ satisfying the properties that $(a,b)\sim_\ell(a',b')$ for some $\ell$ in $\L_2$ and further, there is an annulus $(\gamma,\gamma')$ containing both segments $(a,b)$ and $(a',b')$, where $\gamma,\gamma'\in \Gamma_c$ for some $c$ in $C \setminus C_1$ or $\gamma,\gamma'\in \Gamma_v$ for some $v$ in $V \setminus V_1$ (since $\ell\in\L_2$). 

    We bound the number of such pairs before dealing with the term that arises from Theorem~\ref{thm:pt-plane-restricted}. We restrict our attention to the case where $\gamma, \gamma'$ are circles, the case where they are lines is similar. Summing over $r\not=0$, we show show that the total number of such pairs is $\ll p K |A|$. In fact we count quintuples $(a,b,a',b',c)$ where $c$ is the common center of $\gamma$ and $\gamma'$ (that they share a centre was shown in \cite[Lemma 2.2]{lund2018pinned}). The picture in \cite[Lemma 3.2]{lund2016bis} might be helpful.
    
    The number of pairs $(a,c)$ is at most $2|A|$ because $a \in A_\gamma$ for some $\gamma \in \Gamma$ (and knowing $\gamma$ gives $c$) and by Lemma~\ref{lem:Rich}
    \[
        \sum_{\gamma \in \Gamma} |A_{\gamma}| \leq 2 \sum_{\gamma \in \Gamma} \left| A_{\gamma} \setminus \bigcup_{\gamma'' \neq \gamma} A_{\gamma''} \right| \leq 2 |A|.
    \]
    For each $(a,c)$ there number of admissible $b$ is at most $|A_c| \leq K$. For each $(a,b,c)$ there are at most $(p+1)$ possibilities for $\ell$ and therefore at most $(p+1)$ admissible $(a',b')$. Overall there are at most $2(p+1) K |A|$ quintuples $(a,b,a',b'c)$ and therefore $\ll pK|A|$ pairs $(a,b)$ and $(a',b')$.

    Via \eqref{eq: E} and \eqref{eq: IL} we get
    \begin{align} \label{B2}
      \B_2^*(A) &\ll |A|^2 + \sum_{r \neq 0} |S_r|^{3/2} + |A|^{1/2} \sum_{r \neq 0} |S_r| + p K |A| \nonumber \\
      &\leq |A|^2 + \sum_{r \neq 0} |S_r|^{3/2} + |A|^{5/2} + p K |A| \nonumber \\
      &\leq \sum_{r \neq 0} |S_r|^{3/2} + p K |A| + |A|^{5/2} .
    \end{align}
         
    We bound the first summand like we did in the proof of Theorem~\ref{thm:isosceles F}.
    \begin{align*}
        \sum_{r \neq 0} |S_r|^{3/2}  
        & \leq \left( \sum_{r} |S_r| \right)^{1/2} \left( \sum_{r\neq 0} |S_r|^2 \right)^{1/2}\\
        & \leq |A| \left( |A| \T(A) \right)^{1/2}\\
        & = |A| \left( |A| \left(\T(A) - \frac{|A|^3}{p}\right) + \frac{|A|^4}{p} \right)^{1/2}.
    \end{align*}
    
    By the hypothesis in the statement of the Claim~\ref{T2} we get
    \[
        \sum_{r \neq 0} |S_r|^{3/2} \ll |A|^{3/2} \sqrt{\T_{2,\text{bal}}(A)} + \frac{|A|^3}{p^{1/2}}.
    \]
    
    Substituting into \eqref{B2} yields
    \[
        \B_2^*(A) \ll |A|^{3/2} \T_{2,\text{bal}}(A)^{1/2} + \frac{|A|^3}{p^{1/2}} + |A|^{5/2} + p K |A|.
    \]
    Substituting this into \eqref{T2*} yields
    \[
    \T_{2,\text{bal}}(A) \ll p^{1/2} |A|^{5/4} \T_{2,\text{bal}}(A)^{1/4} + p^{1/4} |A|^2 + p^{1/2} |A|^{7/4} + p K^{1/2} |A|.
    \]
    Standard calculations yield the claimed bound on $\T_{2,\text{bal}}(A)$.

\appendix

\section{Clifford algebra computations}
\label{sec:clifford-calc}
    This section is a short digest of Clifford algebras over fields.
    We follow Klawitter and Hagemann \cite{klawitter2013kinematic}, who give a similar exposition for Clifford algebras over $\R$. 
    A similar argument also applies to the symmetry group of the sphere, see the Appendix in the paper \cite{rudnev2016klein}  by the fourth author and Selig.
        
    For a vector space $V$ with a quadratic form $Q$, the Clifford algebra $\cl(V,Q)$ is the largest algebra containing $V$ and satisfying the relation that $x^2=Q(x)$ for all $x\in V$, where $x^2$ is the square of $x$ in the algebra.
(Recall that an algebra over a field $\F$ is a ring with a homomorphism from $\F$ into its center.)
We use $e_0$ to denote the multiplicative identity of $\cl(V,Q)$; the field $\F$ is embedded in $\cl(V,Q)$ by $x\mapsto x e_0$.

    If $V$ is an $n$-dimensional vector space over a finite field $\F$ of odd characteristic, then there is a basis $e_1,\ldots,e_n$ of $V$ such that $Q(e_i)=\lambda_i$, where $\lambda_i$ is one of: 0, 1, a non-square.
    This basis is orthogonal with respect to the bilinear form associated to $Q$.
    The Clifford algebra is a $2^n$-dimensional $\F$ vector space with basis $e_{i_1\ldots i_k}$ where $i_1<\cdots<i_k$ and $0\leq k\leq n$ defined by
    \begin{equation}
      \label{eq:clifford-basis}
              e_{i_1\ldots i_k}=e_{i_1}\cdots e_{i_k}.
    \end{equation}
The rules for multiplication in $\cl(V,Q)$ are given by $e_ie_j=-e_je_i$ for $i\not=j$ and $e_i^2=\lambda_i$, where $1\leq i,j\leq n$; these rules extend to all of $\cl(V,Q)$ by \eqref{eq:clifford-basis} and linearity.
    
The Clifford algebra $\cl(V,Q)$ splits as a direct sum of exterior products
    \[  
        \cl(V,Q)=\bigoplus_{i=0}^n \bigwedge^i V
    \]
and is $\Z/2$-graded:
    \[
        \cl(V,Q)=\cl(V,Q)^+\oplus \cl(V,Q)^-,
    \]
where
    \[
        \cl(V,Q)^+:=\bigoplus_{\substack{i=0\\ i\equiv 0\pmod 2}}^n \bigwedge^i V
            \andd
        \cl(V,Q)^-:=\bigoplus_{\substack{i=0\\ i\equiv 1\pmod 2}}^n \bigwedge^i V.
    \]
    The dimension of the even subalgebra $\cl(V,Q)^+$ is $2^{n-1}$.
We identify $\bigwedge^0 V$ with $\F$ and $\bigwedge^1 V$ with $V$.
    
We define two involution of $\cl(V,Q)$.
The first, called \emph{conjugation}, is denoted by an asterisk.
For the basis elements of $V$ we define conjugation by $e_i^*=-e_i$.
We extend conjugation to other basis elements by changing the order of multiplication
    \[
        (e_{i_1}e_{i_2}\cdots e_{i_k})^*:=(-1)^ke_{i_k}\cdots e_{i_2}e_{i_1}\quad 0\leq i_1<i_2<\cdots<i_k\leq n.
    \]
Finally, we extend conjugation to $\cl(V,Q)$ by linearity.
One can check that $(\fa\fb)^*=\fb^*\fa^*$ for an elements $\fa,\fb\in\cl(V,Q)$.
(Notice that if $\fa\in\bigwedge^k V$, then $\fa^*=(-1)^{k(k+1)/2}\fa$.)
For any element $\fa$ in $\cl(V,Q)$, the product $\fa\fa^*$ is a scalar.
We define the \emph{norm}\/ of an element $\fa$ by $N(\fa)=\fa\fa^*$; notice that $N(\fa\fb)=N(\fa)N(\fb)$.
    
    The second involution of $\cl(V,Q)$, called the \emph{main involution}, is denoted by $\alpha$ and defined by $\alpha(e_i)=-e_i$ and extended to $\cl(V,Q)$ by linearity and the rules for multiplication.
The main involution is an algebra automorphism.
    Clearly $\alpha$ fixes the even subalgebra $\cl(V,Q)^+$ and acts by multiplication by $-1$ on the odd subalgebra $\cl(V,Q)^-$.

    Let $\cl^\times(V,Q)$ denote the set of invertible elements of $\cl(V,Q)$, which we call \emph{units}.
    If $\fa\in V$ and $N(\fa)\not=0$, then $\fa\in\cl^\times(V,Q)$, and $\fa^{-1}=\fa^*/N(\fa)$.
    The \emph{Clifford group}\/ associated to $\cl(V,Q)$ is defined by
    \[
        \Gamma(\cl(V,Q)):=\left\{\fg\in\cl^\times(V,Q)\colon \alpha(\fg)V\fg^{-1}\subseteq V \right\}.
    \]
    We say that the map $\fv\mapsto \alpha(\fg)\fv\fg^{-1}$ is the \emph{sandwich operator}\/ associated to an element $\fg\in\Gamma(\cl(V,Q))$.
    
\newcommand{\SOQ}{\ensuremath{\mathrm{SO}(Q_0)}}
\newcommand{\SFQ}{\ensuremath{\mathrm{SF}(Q_0)}}
    Given a quadratic form $Q_0$ on $\F^2$ with $Q_0(e_1)=1$ and $Q_0(e_2)=-\lambda$, let $\SOQ$ denote the group of rotations preserving $Q_0$:
    \[
        \SOQ:= \left\{
        \begin{pmatrix}
            u & v \\
            \lambda v & u\\
        \end{pmatrix}
        \colon u^2-\lambda v^2=1
        \right\},
    \]
    and let $\SFQ$ denote the group of rigid motions of $\F^2$ generated by $\SOQ$ and the group of translations.
If $\lambda=-1$, then $\SOQ=\SO$ and $\SFQ=\SF$.
The elements of $\SFQ$ can be represented by matrices of the form
\begin{equation}
  \label{eq:SFQ-matrix-form}
  \begin{pmatrix}
    u & v & s \\
    \lambda v & u & t\\
    0 & 0 & 1\\
  \end{pmatrix}.
\end{equation}
The action of an element of $\SFQ$ on a vector $(x,y)^T$ in $\F^2$ corresponds to matrix multiplication on the vector $(x,y,1)^T$.

    \begin{prop}
      \label{prop:1}
        Let $V=\F^3$ and define $Q$ on $V$ by $Q(x,y,z)=Q_0(x,y)$, let $G=(\cl(V,Q)^+)^\times$ be the group of units of the even subalgebra, and let $Z$ be its centre.
        Then $G/Z$ is isomorphic to $\SFQ$.
    \end{prop}
The main idea in this proof is that the set of units in the even subalgebra of $\cl(V,Q)$ act on $V$ by sandwich operator; this gives us a matrix representation of the group of units, and the dual of the representation is precisely the matrix representation of $\SFQ$.
    \begin{proof}
        By our definition of $Q$, $e_1^2=1,e_2^2=-\lambda,e_3^2=0$, $\cl(V,Q)$ is spanned by
        \[
            e_0, e_1,e_2,e_3,e_{12},e_{13},e_{23},e_{123},
        \]
        and $\cl(V,Q)^+$ is spanned by $e_0,e_{12},e_{13},e_{23}$.
        If $\fg= g_0e_0 + g_{12}e_{12}+g_{13}e_{13}+g_{23}e_{23}$, then
        \[
            N(\fg)=\fg\fg^*=g_0^2-\lambda g_{12}^2\,.
        \]
        Thus, if $g_0^2-\lambda g_{12}^2\not=0$, the inverse of $\fg$ is 
        \[
            \fg^{-1}=\frac 1{g_0^2-\lambda g_{12}^2}\fg^*.
        \]
        This determines the group of units explicitly.
    
        One can show by a computation that $G$ acts on $V$ by the sandwich product $(\fg,\fv)\mapsto \fg\fv\fg^{-1}$ (that is $G=\Gamma(\cl(V,Q)^+)$).
        In fact, the action of general element $\fg=g_0e_0 +g_{12}e_{12}+g_{13}e_{13}+g_{23}e_{23}$ in $G$ is given by
        \begin{align*}
            \fg e_1 \fg^{-1}&=\frac{g_0^2+\lambda g_{12}^2}{g_0^2-\lambda g_{12}^2}e_1
            +\frac{-2g_0g_{12}}{g_0^2-\lambda g_{12}^2}e_2
            +\frac{-2(g_0g_{13}+\lambda g_{12}g_{23})}{g_0^2-\lambda g_{12}^2}e_3,\\
            \fg e_2 \fg^{-1}&=\frac{-2\lambda g_0g_{12}}{g_0^2-\lambda g_{12}^2}e_1
            +\frac{g_0^2+\lambda g_{12}^2}{g_0^2-\lambda g_{12}^2}e_2
            +\frac{2\lambda(g_0g_{23}+ g_{12}g_{13})}{g_0^2-\lambda g_{12}^2}e_3,\\
            \fg e_3 \fg^{-1}&=e_3.
        \end{align*}
        Let $\rho\colon G\to\GL(V)$ denote this representation.

        The \emph{dual representation}\/ $\rho^*(\fg):=\rho(\fg^{-1})^T$, where $T$ denotes the transpose, acts on the dual space $V^*$, and in the standard basis $\{f_1,f_2,f_3\}$ on $V^*$ defined by $f_i(e_j)=\delta_{ij}$, we have
        \begin{equation}
          \label{eq:contragredient}
            \rho^*(\fg^{-1})=
            \frac 1{g_0^2-\lambda g_{12}^2}
            \begin{pmatrix}
                g_0^2+\lambda g_{12}^2 &-2g_0g_{12} & -2(g_0g_{13}+\lambda g_{12}g_{23})\\
                -2\lambda g_0g_{12} & g_0^2+\lambda g_{12}^2 & 2\lambda(g_0g_{23}+ g_{12}g_{13})\\
                0 & 0&  g_0^2-\lambda g_{12}^2\\
            \end{pmatrix}.
        \end{equation}
        The kernels of $\rho$ and $\rho^*$ are both equal to the subgroup $Z:=\{g_0e_0\colon g_0\not=0\}$.

 We wish to show that $G/Z$ is isomorphic to $\SFQ$.
By inspection, the image of $G$ under $\rho^*$ consists of matrices of the form~\eqref{eq:SFQ-matrix-form}, so $G/Z$ is isomorphic to a subgroup of $\SFQ$, so it remains to show that $\rho^*$ is surjective.

Let $R$ be the subgroup defined by $g_{13}=g_{23}=0$; the rational parameterisation of the ellipse defined by $u^2-\lambda v^2=1$ shows that $\rho^*(R)$ maps onto the subgroup $\SOQ\subseteq\SFQ$.
(By rational parameterisation, we mean the map
\[
t\mapsto \left( \frac{t^2+\lambda}{t^2-\lambda},\frac{-2t}{t^2-\lambda} \right),
\]
which is a bijection from $\F$ onto $\{(u,v)\colon u^2-\lambda v^2=1\}\setminus\{(1,0)\}$.
A similar argument is used by Bennett, Iosevich, and Pakianathan \cite{bennett2014three-point}.)

On the other hand, it is clear that the subgroup $T$ defined by $g_0=1,g_{12}=0$ is bijective with the translation subgroup of $\SFQ$.
Since these subgroups generate $\SFQ$, we see that $\rho^*$ is surjective.
\end{proof}
    
    We have shown more: $\SFQ$ is naturally identified with a (Zariski open) subset of $\PF$, and the nature of this identification yields some desirable features.
In particular, as in Corollary~\ref{cor:transporters}, the set of transformations in $\SFQ$ that map a point $x\in \F^2$ to a point $y\in\F^2$ is a projective line. 
    
    Let $\kappa\colon\SFQ\to G/Z$ denote the inverse of $\rho^*\colon G/Z\to\SFQ$.
    This is the \emph{kinematic mapping}\/ of Blaschke and Gr\"unwald, who both sought to embed the group of rigid motions in projective space.
    Let $\PF$ denote projective three space; we write $[X_0\colon X_1\colon X_2\colon X_3]$ for a typical point of $\PF$.
    \begin{cor}
      \label{cor:1}
        There is a bijection $\kappa\colon\SFQ\to\PF\setminus\{X_0^2-\lambda X_1^2=0\}$ such that the image of the rotation subgroup and translation subgroups are projective lines.
    
        Further, for all $g\in\SFQ$ there are projective maps $\phi_g\colon \PF\to \PF$ and $\phi^g\colon \PF\to \PF$ such that for all $x\in\SFQ$
        \[
            \kappa(gx) =\phi_g(\kappa(x))\andd \kappa(xg)=\phi^g(\kappa(x)).
        \]
            \end{cor}
    \begin{proof}
        The even subalgebra $\cl(V,Q)^+$ is isomorphic to $\F^4$ as a vector space, so the projective space $\PP(\cl(V,Q)^+)$ is $\PF$.
        On the other hand, $\PP(\cl(V,Q)^+)$ is just $\cl(V,Q)^+$ modulo the action of the multiplicative subgroup $Z$, so we have $G/Z\subseteq\PP(\cl(V,Q)^+)$.
        In fact, $G/Z$ consists of all points $[g_0\colon g_{12}\colon g_{13}\colon g_{23}]$ such that $g_0^2-\lambda g_{12}^2\not=0$.
        
        Since $\cl(V,Q)^+$ is an $\F$-algebra, left and right multiplication are $\F$-linear transformations.
        That is, if $\tilde\phi_{\fg}(\fv):=\fg\fv$ and $\tilde\phi^{\fg}(\fv):=\fv\fg$, then $\tilde\phi_{\fg}$ and $\tilde\phi^{\fg}$ are linear transformations.
        It follows that left and right translation in $G/Z$ are \emph{projective transformations}\/ of $\PF$.
        
    \end{proof}

\paragraph{\textbf{Acknowledgments}.} The first and fourth listed authors are partially supported by the Leverhulme Trust Grant RPG--2017--371. Parts of the project begun while they visited The University of Georgia. The second listed author is supported by the NSF Award 1723016 and gratefully acknowledges the support from the RTG in Algebraic Geometry, Algebra, and Number Theory at the University of Georgia, and from the NSF RTG grant DMS-1344994. The third listed author is supported by Swiss National Science Foundation grants P400P2-183916 and P4P4P2-191067. The fifth listed author is supported by Austrian Science Fund FWF Project P 30405-N32.

We thank Brandon Hanson, Alex Iosevich, David Ellis and Konrad Swanepoel for comments we received on various stages of writing this paper, and extend our thanks to the anonymous referee for a thorough report and helpful suggestions. Last but not least, we owe much credit to Jon Selig, who once introduced the third author to the concepts of geometric kinematics.
\bibliographystyle{siam}
\bibliography{library}
\end{document}